\documentclass[12pt]{article}
\usepackage{amsmath,bbm,array,amsfonts,float,mathtools,amsthm}
\usepackage{jheppub}
\usepackage{caption,subcaption,multirow}
\usepackage[toc,page]{appendix}
\usepackage[utf8]{inputenc}
\usepackage[autostyle=true]{csquotes}

\numberwithin{equation}{section}
\numberwithin{figure}{section}
\numberwithin{table}{section}
\usepackage{todonotes}

\title{Cluster Algebras: Network Science and Machine Learning}

\author[a,b,c]{Pierre-Philippe Dechant,}
\author[d,e,f,g]{Yang-Hui He,}
\author[e,d]{Elli Heyes,}
\author[e,d]{Edward Hirst}

\affiliation[a]{
	School of Science, Health and Technology, York St John University, YO31 7EX, UK}
\affiliation[b]{
	Department of Mathematics, University of York, YO10 5DD, UK}
\affiliation[c]{
	York Cross-disciplinary Centre for Systems Analysis, University of York, YO10 5DD, UK}
\affiliation[d]{
    London Institute for Mathematical Sciences, Royal Institution, London W1S 4BS, UK}
\affiliation[e]{
	Department of Mathematics, City, University of London, EC1V 0HB, UK}
\affiliation[f]{
	Merton College, University of Oxford, OX1 4JD, UK}
\affiliation[g]{
	School of Physics, NanKai University, Tianjin, 300071, P.R. China}

\emailAdd{ppd22@cantab.net}
\emailAdd{hey@maths.ox.ac.uk}
\emailAdd{elli.heyes@city.ac.uk}
\emailAdd{edward.hirst@city.ac.uk}

\preprint{
	\begin{flushright}
		LIMS-2022-011
	\end{flushright}
}
\abstract{
Cluster algebras have recently become an important player in mathematics and physics.
In this work, we investigate them through the lens of modern data science, specifically with techniques from network science and machine learning.
Network analysis methods are applied to the exchange graphs for cluster algebras of varying mutation types. 
The analysis indicates that when the graphs are represented without identifying by permutation equivalence between clusters an elegant symmetry emerges in the quiver exchange graph embedding. 
The ratio between number of seeds and number of quivers associated to this symmetry is computed for finite Dynkin type algebras up to rank 5, and conjectured for higher ranks. 
Simple machine learning techniques successfully learn to classify cluster algebras using the data of seeds. 
The learning performance exceeds 0.9 accuracies between algebras of the same mutation type and between types, as well as relative to artificially generated data.\\

\emph{To the memory of Professor John K. S. McKay (1939-2022), with deepest respect.}
}

\begin{document}
\maketitle

\section{Introduction}
Cluster algebras \cite{CA_1,CA_2,williams2013cluster,marsh2013lecture,cheung2021cluster,Duan_2020}, have seen an exponential growth in interest over recent years. 
As objects which take root in combinatorics, geometry, and number theory, their mutation structure connecting sets of algebra generators is arising in more and more new contexts.
Of these contexts, physics has offered a significant number. 
Beyond an interpretation of the exchange graph mutation process as the action of Seiberg duality \cite{Seiberg_1995} in connecting IR equivalent 4d $\mathcal{N}=1$ gauge theories \cite{Feng:2001bn,Fordy:2009qz,Benini:2014mia,Franco:2014nca,Franco:2003ja}, there are even extensions of these algebras to incorporate more general quantum field theory symmetries \cite{Franco:2017lpa}.
Cluster algebras have led to a variety of results in higher Toda theories \cite{Fock2003ModuliSO,Fock2009QD,Williams:2014efa}, and provided new methods of looking at wall crossing \cite{Gaiotto:2011tf,Kontsevich:2008fj}.
Moreover, a particularly fruitful application has been to the computation of scattering amplitudes \cite{Golden:2013xva,Golden_2014,He2021CA}.
Furthermore, the perturbative and computationally demanding approach of computing Feynman integrals order by order can in some cases be replaced through consideration of amplituhedra and positive Grassmannians \cite{Arkani-Hamed:2013jha,Arkani-Hamed:2012zlh}.

As the theory develops and computations require knowledge of larger and larger cluster algebras, the field of their physical application becomes better and better suited for `big data' techniques.
Machine learning (ML) is a topical and equally fast-developing new field of techniques suited for learning properties of large datasets, especially in theoretical physics and pure mathematics.
The techniques for ML range from non-linear function fitting and landscape searching, to clustering and pattern recognition; and in \cite{He:2017aed,He:2017set,Carifio:2017bov,Krefl:2017yox,Ruehle:2017mzq} they were first introduced to the string landscape.
Further to their success in the algebraic geometry sector of string theory \cite{Bull:2018uow,He:2020lbz,Jejjala:2020wcc,Anderson:2020hux,ashmore2021machine,Cole:2021nnt,gao2021machine,Berglund:2021ztg,Berman:2021mcw,He:2018jtw,Abel:2021ddu,Jejjala:2022lxh} and the related high-energy physics \cite{Krippendorf:2020gny,Halverson:2021aot,Bao:2021auj,Bao:2021olg,Arias-Tamargo:2022qgb,Hirst:2022qqr}, strong results from ML application have also been seen in various fields of mathematics \cite{He:2018jtw,Jejjala:2019kio,Gukov:2020qaj,He:2020fdg,He:2020eva,He:2021oav,Bao:2021ofk,davies2021advancing,DLGMC,Bena:2021wyr}.
Particularly relevant to the present context are
\cite{He:2019nzx} where the study of ML on algebraic structures was initiated, as well as \cite{Amoros:2021rrx} where ML was utilized in classification problems in commutative algebra and \cite{osti_10276457} where learning strategies were imposed in the key step of Buchberger algorithm in algebraic geometry.

In this work, we initiate the application of ML techniques to the theory of cluster algebras, with the hope that as they develop they will offer insight beyond the successes here.
Importantly we extend the work in \cite{Bao:2020nbi}, where ML was applied to quiver mutation, to now examine the full cluster algebra seeds.
This complements further extension of the databases considered to more general skew-symmetrisable quivers, as well as broadening the range of analysis techniques to include network science also.
We begin in section \S\ref{CAreview} with a brief review of cluster algebras, the relevant terminology, and how they naturally fall into classification types. 
A key point of focus is the exchange graph, which dictates the structure of the mutation action, and in section \S\ref{dataanalysis} we examine these for a variety of cluster algebras using a range of techniques from the network analysis toolbox.
Remarkably, it is discovered that by \textit{not} identifying cluster seeds which are equivalent under permutation, the finite type cluster algebras have ratios between the total number of distinct seeds and distinct quivers which are strictly integers.
In section \S\ref{ML} the prototypical supervised machine learning architecture of feed-forward neural networks is used to learn structures inherent in seeds and clusters. Binary as well as multi-class classifications between the different algebras considered show very strong performance.
Summary comments and outlook are provided in section \S\ref{summary}.

The coding scripts and datasets associated to this work are available at the respective GitHub: \url{https://github.com/edhirst/ClusterAlgebrasML.git}.

\section{Cluster Algebras Review}\label{CAreview}
A \textit{Cluster Algebra} of rank $r$ is a commutative subring of an ambient field of rational functions in $r$ variables \cite{CA_1,CA_2,williams2013cluster}. 
Each algebra is constructively defined through some size $r$ subset of the algebra's generators $\{x_i\}$ known as a \textit{cluster} (where each generator is a \textit{cluster variable}), along with an $r \times r$ \textit{exchange matrix}. 
Together, the cluster and the exchange matrix form a \textit{seed}, from which all cluster variables can be generated via \textit{mutation}, thus providing the full set of generators to define the algebra.

From an initial seed, action of the \textit{mutation} process on a chosen cluster variable in a cluster exchanges it for another cluster variable not contained in that cluster. 
The exchange matrix also gets mutated.
Mutation hence changes the seed to another one, 
different from the original. 
Continuing this, repeated action of the mutation process would produce all the seeds (of which there may be infinite) and the union of all the clusters from all the seeds is the algebra's complete set of generators.
Hence any seed can be validly used to generate all other seeds via mutation.

The exchange matrix associated to each seed may be defined via the signed adjacency matrix of a \textit{quiver}, making it a skew-symmetric matrix, $b_{ij}$. Note that this does restrict the cluster algebra definition to skew-symmetric cluster algebras of geometric type. These quivers are weighted directed graphs on $r$ vertices with no loops or 2-cycles such that $\ell$ arrows exist from vertex $i$ to vertex $j$ if $b_{ij} = -b_{ji} = \ell$ in the exchange matrix.
An extension can be made to skew-symmetrizable but not skew-symmetric matrices, drawing the respective quiver edges with a double weighting $(b_{ij},\ b_{ji})$ where $b_{ij}\neq -b_{ji}$.

Given a seed, one can mutate about each of the $r$ cluster variables, each choice producing a different seed. 
For mutation about variable $x_k$ the cluster $\{x_i\}$ updates its variables via
\begin{equation}\label{eq:clustermutation}
    x_i \longmapsto x'_i = 
        \begin{cases}
        \big(\prod_{b_{\mu i}>0}x_\mu^{b_{\mu i}} + \prod_{b_{\nu i}<0}x_\nu^{-b_{\nu i}} \big) / x_i \quad & \big| \ i=k\\
        x_i & \big| \ i \neq k
        \end{cases}
\end{equation}
such that only the chosen cluster variable updates, and it does so based on which variables' vertices are connected by inward ($b_{\mu k}>0$) or outward ($b_{\nu k}<0$) edges to the mutated vertex in the exchange matrix's quiver.
This process produces a new cluster variable which is a Laurent polynomial in the initial seed's variables, and this Laurent phenomenon is preserved such that all cluster variables after any number of mutations can be expressed as Laurent polynomials of these initial cluster variables.

Further to changing the cluster, mutation of a seed about a chosen cluster variable, $x_k$, also mutates the seed's exchange matrix
\begin{equation}\label{eq:quivermutation}
    b_{ij} \longmapsto b'_{ij} =
        \begin{cases}
        \quad -b_{ij} & \big| \ k = i \ \text{or} \ j\\
        b_{ij} + b_{ik}b_{kj} & \big| \ b_{ik} > 0 \ \text{and} \ b_{kj} > 0\\
        b_{ij} - b_{ik}b_{kj} & \big| \ b_{ik} < 0 \ \text{and} \ b_{kj} < 0\\
        \quad \ \ b_{ij} & \big| \ \text{otherwise}
        \end{cases}
\end{equation}
This in the language of the quiver amounts to reversing the direction of all arrows incident to vertex $k$ associated to variable $x_k$, then for all 2-paths through the vertex $\mu \rightarrow k \rightarrow \nu$ adding an arrow $\nu \rightarrow \mu$ (closing the 2-path into a 3-cycle), and finally removing all arrows that form 2-cycles in the quiver.

In the physical interpretation, quivers may be used to represent the particle content of supersymmetric gauge theories \cite{Douglas:1996sw,Bao:2021ohf}.
Here quiver vertices dictate the gauge groups of the theory (which has 4 supercharges when a directed graph) along with the respective vector multiplets, whilst arrows represent chiral multiplets in bifundamental representations of the gauge groups from the vertices that they connect. 
The mutation process on the quivers then is equivalent to Seiberg duality \cite{Seiberg_1995,Feng_2001,Franco:2003ja,Franco_2004,Alim2013BPS}, such that the gauge theories of all quivers that are related via mutation have the same IR fixed point.
The reader is referred to \cite{Franco:2003ja,Franco_2004} for the duality trees in gauge theory.

\subsection{Exchange Graphs}
Since an initial seed can be mutated on each of its $r$ cluster variables, each seed can mutate to $r$ different seeds. 
Iterating this mutation produces all possible seeds; each may be used to define the entire cluster algebra (hence the name). 
This mutation process between seeds may also be represented in the language of graph theory. 
Here each seed is associated to a vertex, and an edge exists between any two vertices if mutation takes either of the vertices' seeds to the other vertex's seed.
Since the mutation process is an involution (such that mutating on a variable just produced via mutation reverts it back to the original variable) all the edges are undirected.

The \textit{exchange graph} of a cluster algebra is the graph representing the relationships between all the seeds, where connectivity is given by mutation. It is generated by inputting an initial seed, mutating this on all $r$ of its variables to give all the `depth 1' seeds (i.e. 1 mutation away from the initial seed), then mutating each of these new seeds on each of their $r-1$ variables that weren't mutated on from the previous depth to give the depth 2 seeds, then iterating this process either infinitely or until no new seeds are generated.
Note that the exchange graph is \textit{not} the graph formed from an exchange matrix; this would be a quiver. 

One may also construct another graph representing connectivity under mutation by considering only the quiver part of the seeds, i.e. the exchange matrix. This exchange graph represents how all the possible exchange matrices are connected via mutation action. This in turn is called the \textit{quiver exchange graph}, and is often significantly smaller than the exchange graph with cluster information at any depth. This is because many seeds may have the same exchange matrix but different clusters, making them the same vertex in the quiver exchange graph but different vertices in the exchange graph.
To make clear the differentiation between these types of exchange graph we will call the graph where the seeds include the cluster information the \textit{seed exchange graph}.

In both cases, as the seed exchange graph / quiver exchange graph is generated, seeds / quivers produced at the next depth may have already occurred at previous depths. 
These vertices are then combined; that way cycles can be formed in the exchange graph.
The frequency of encountering previous seeds strongly depends on the cluster algebra type, and analysis of the occurrence of these cycles for all types will be a central theme to this work.

\begin{figure}[!tb]
	\centering
	\begin{subfigure}{0.45\textwidth}
    	\centering
    	\includegraphics[width=0.9\textwidth]{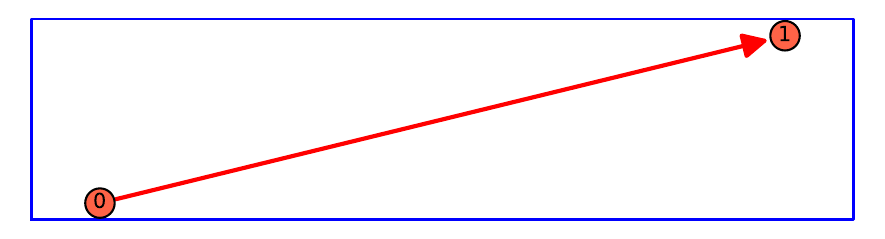}
    	\caption{A2 quiver}\label{fig:A2quiv}
	\end{subfigure} 
    \begin{subfigure}{0.45\textwidth}
    	\centering
    	\includegraphics[width=0.65\textwidth]{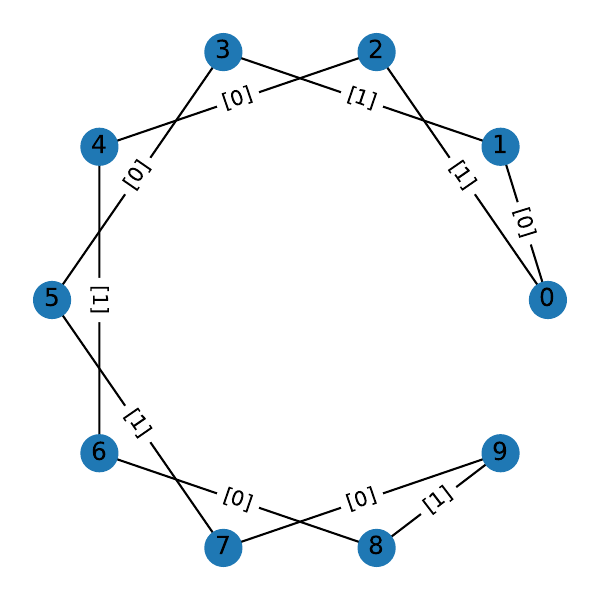}
    	\caption{A2 seed exchange graph}\label{fig:A2EG}
    \end{subfigure} \\
\caption{The quiver for the $A_2$ example cluster algebra (a), as well as its exchange graph (b) where permutation equivalence is not applied. The respective clusters for each vertex in the exchange graph are: $\{ 0: [x_1,x_2], \ 1: [(x_2 + 1)/x_1, x_2], \ 2: [x_1, (x_1 + 1)/x_2], \ 3: [(x_2 + 1)/x_1, (x_1 + x_2 + 1)/(x_1x_2)], \ 4: [(x_1 + x_2 + 1)/(x_1x_2), (x_1 + 1)/x_2], \ 5: [(x_1 + 1)/x_2, (x_1 + x_2 + 1)/(x_1x_2)], \ 6: [(x_1 + x_2 + 1)/(x_1x_2), (x_2 + 1)/x_1], \ 7: [(x_1 + 1)/x_2, x_1], \ 8: [x_2, (x_2 + 1)/x_1, \ 9: [x_2, x_1] \}$, and the vertex mutated on to connect each seed is given as the respective edge feature. Note that throughout this paper vertices are labelled in order of generation, here this $A_2$ seed exchange graph is drawn with vertices in order of labelling, but may be be unfolded to show its single cycle nature with a different vertex ordering.}\label{fig:A2example}
\end{figure}

\paragraph{Example: $A_2$}
The prototypical example for cluster algebras is the finite-type rank 2 algebra: $A_2$.
This has quiver given by the directed $A_2$ Dynkin diagram shown in Figure \ref{fig:A2quiv}.

Starting with the initial seed $\{x_1,x_2\}$, one can mutate about either of the vertices.
With each mutation, since the rank is too small for the quiver to include 2-paths, the only change to the quiver is a reversal of the single edge's orientation.
Therefore in this simple $A_2$ case the quiver exchange graph is just 2 vertices connected by an edge to represent both these quivers (which we consider not equivalent via permutation), the second quiver exchange matrix being just the transpose of the first. 
Furthermore in the seed exchange graph, Figure \ref{fig:A2EG}, the quiver alternates between these two forms around the exchange graph's loop.

To exemplify the mutation process, consider for this $A_2$ algebra the mutation of the seed associated to vertex 4 in the exchange graph on the quiver vertex labelled 1 (associated to the second variable), which mutates the seed to the seed labelled 6 in the seed exchange graph.
The quiver associated to this seed is the same as for the initial seed, as shown in Figure \ref{fig:A2quiv}. 
Therefore the exchange matrix is \footnotesize $\begin{pmatrix}0 & 1\\ -1 & 0 \end{pmatrix}$\normalsize, which updates to its transpose under the process in equation \eqref{eq:quivermutation}.
The cluster $[(x_1 + x_2 + 1)/(x_1x_2), (x_1 + 1)/x_2]$ is respectively mutated on its second variable following the process in equation \eqref{eq:clustermutation} with $k=2$. 
The process hence keeps the first variable unchanged and updates the second variable where the first numerator term in \eqref{eq:clustermutation} is just the first variable (as $EM_{12}=1$) and second numerator term is just 1 (since there are no negative entries in the second column).
Therefore $\frac{(x_1 + 1)}{x_2} \longmapsto \frac{(x_1 + x_2 + 1)/(x_1x_2)+1}{(x_1 + 1)/x_2} = \frac{(x_2+1)}{x_1}$, matching the expected cluster for seed 6.

\subsection{Cluster Algebra Types}\label{sec:types}
Cluster algebras can be classified into 3 distinct types:\\
$\quad$ 1) \textit{Finite type}: These are defined by having finite numbers of seeds and thus finitely many clusters, exchange matrices, and also cluster variables. All cluster algebras of this type are formed from quivers which take the form of oriented ADE Dynkin diagrams \cite{CA_2,Gabriel}.
Since there are finitely many seeds and exchange matrices both the seed exchange graph and the quiver exchange graph form finite compact polytopes. These polytopes are called generalised associahedra (for this finite case).
This algebra type may also be identified by the condition that $|b_{ij}b_{ji}|\leq 3 \ \forall \ i,j$ in all exchange matrices in all cluster algebra seeds. \\
$\quad$ 2) \textit{Finite-mutation type}: Since the quiver exchange graph is strictly smaller than the seed exchange graph, it can be finite even if the seed exchange graph is infinite. 
Therefore finite-mutation type cluster algebras are defined by having finitely many exchange matrices but may have infinitely many cluster variables, and hence also infinitely many clusters/seeds (i.e. there may be infinitely many seeds with the same exchange matrix but with different cluster variables). 
Therefore whilst the seed exchange graph may be infinite, the quiver exchange graph is finite. These cluster algebras thus naturally include finite type but also others; in particular, the classification consists of cluster algebras formed from rank 2 quivers, from triangulations of marked surfaces, and a set of exceptional cases \cite{Felikson_2012,DerkOwen}.
Extending the finite type condition, finite-mutation type may also be identified as rank 2 or by $|b_{ij}b_{ji}|\leq 4 \ \forall \ i,j$ in all exchange matrices in all the seeds \cite{Felikson_2010}. \\ 
$\quad$ 3) \textit{Infinite type}: These encompass the remaining cases, with infinitely many clusters and exchange matrices leading to infinitely many cluster variables.
Therefore both the seed exchange graph and quiver exchange graph are infinite and non-compact.

As can be seen from the exchange matrix condition, our example $A_2$ cluster algebra is of finite type as there are only 2 exchange matrices both with $|b_{ij}|\leq1 \ \forall \ i,j \ \implies |b_{ij}b_{ji}|\leq 3$.

\section{Exchange Graph Data}\label{dataanalysis}

In this work, we focus on the seed exchange graphs of cluster algebras,
considering a selection that spans the possible types discussed in \S\ref{sec:types}.
We analyse how these graphs take shape as they are generated and introduce the application of ML techniques to study their respective cluster algebras.
In spirit this extends the work in \cite{Bao:2020nbi} where ML was applied to quiver exchange graphs to learn the underlying Seiberg duality.
In order to perform the cluster mutation, and keep track of the seeds computationally, the \texttt{sage} `Cluster Algebra and Quiver' package was used \cite{sage,musiker2011compendium}; the exchange graphs were then represented and analysed with use of the \texttt{python} package \texttt{networkx} \cite{Hagberg2008ExploringNS}; and finally completion of the ML investigations made use of the \texttt{scikit-learn} library \cite{scikit-learn}.

\begin{figure}[H]
	\centering
	\begin{subfigure}{0.45\textwidth}
    	\centering
    	\includegraphics[width=0.8\textwidth]{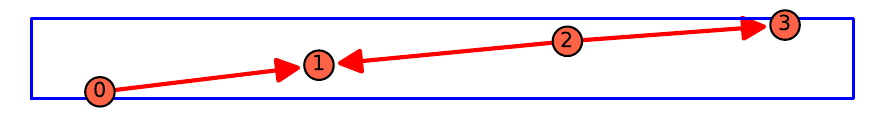}
    	\caption{A4 (finite)}\label{fig:A4}
	\end{subfigure} 
    \begin{subfigure}{0.45\textwidth}
    	\centering
    	\includegraphics[width=0.7\textwidth]{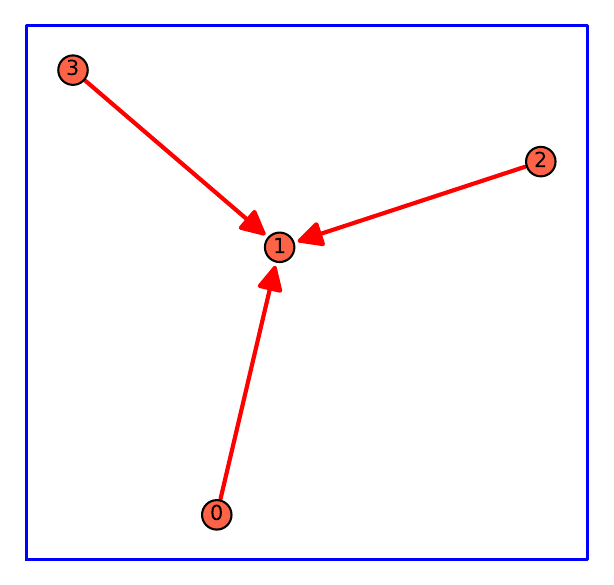}
    	\caption{D4 (finite)}\label{fig:D4}
    \end{subfigure} \\
    \begin{subfigure}{0.45\textwidth}
    	\centering
    	\includegraphics[width=0.8\textwidth]{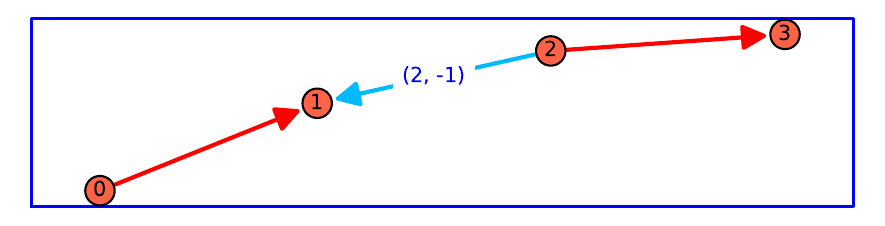}
    	\caption{F4 (finite)}\label{fig:F4}
    \end{subfigure} \\
	\begin{subfigure}{0.45\textwidth}
    	\centering
    	\includegraphics[width=0.7\textwidth]{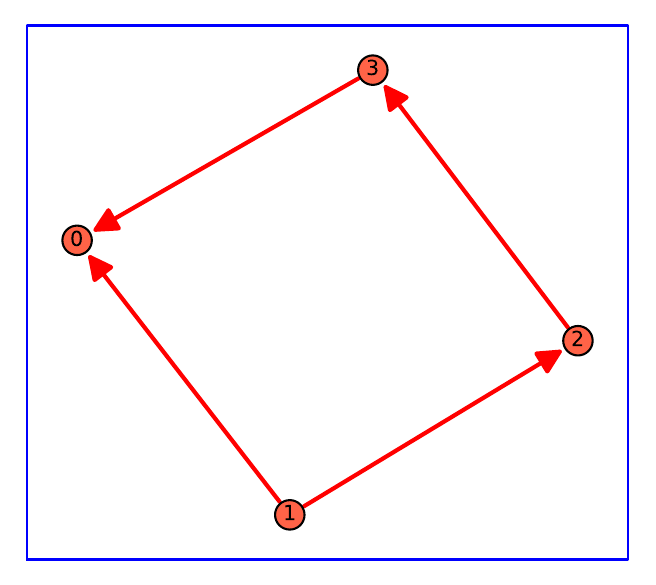}
    	\caption{A13 (finite-mutation)}\label{fig:A13}
	\end{subfigure} 
    \begin{subfigure}{0.45\textwidth}
    	\centering
    	\includegraphics[width=0.6\textwidth]{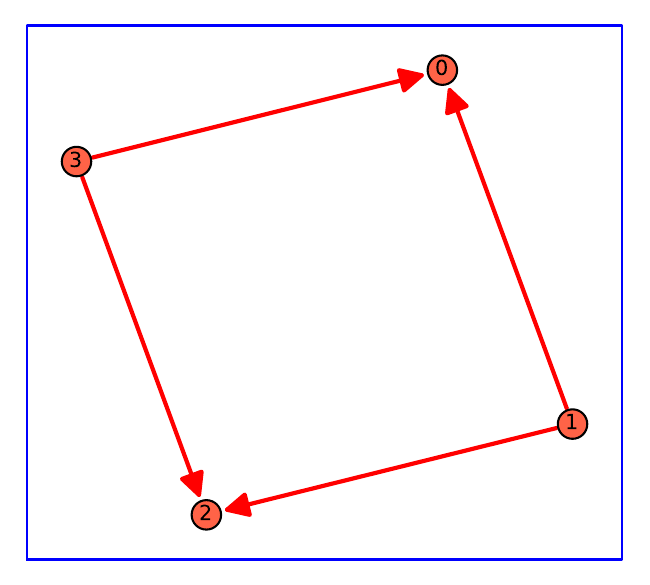}
    	\caption{A22 (finite-mutation)}\label{fig:A22}
    \end{subfigure} \\
	\begin{subfigure}{0.45\textwidth}
    	\centering
    	\includegraphics[width=0.7\textwidth]{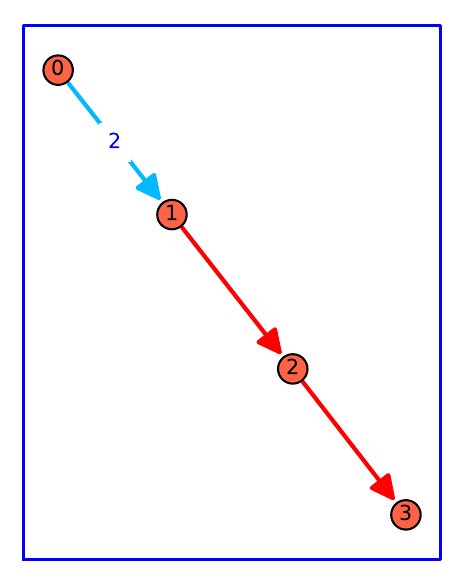}
    	\caption{$\mathcal{I}1$ (infinite)}\label{fig:inf1}
	\end{subfigure} 
    \begin{subfigure}{0.45\textwidth}
    	\centering
    	\includegraphics[width=0.8\textwidth]{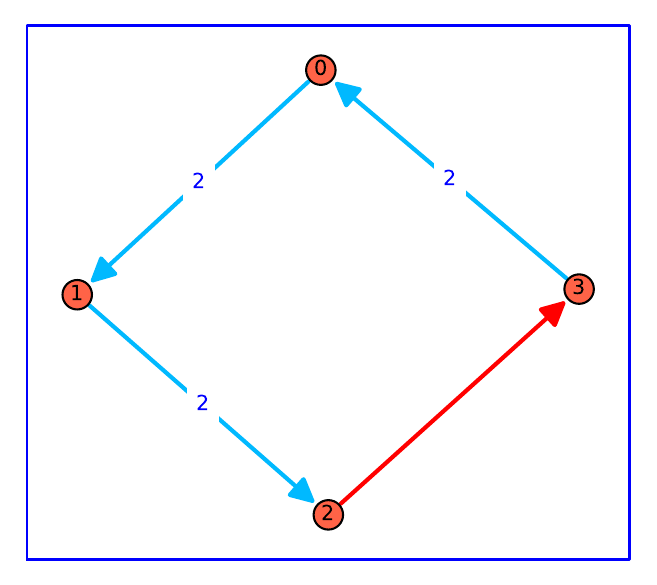}
    	\caption{$\mathcal{I}2$ (infinite)}\label{fig:inf2}
    \end{subfigure} \\
\caption{Quivers defining the exchange matrices for the initial seeds. They are all rank 4, and generate cluster algebras of finite type (a), (b), (c); finite-mutation type that are not finite type (d), (e); and infinite type (f), (g). Vertices are labelled with respect to the row/column number in the exchange matrix; the double edge multiplicity in F4 indicates it is not skew-symmetric.}\label{fig:quivers}
\end{figure}

The cluster algebras considered in these investigations are denoted by the respective initial seeds used to generate them. 
These are hence defined by a choice of exchange matrix, each paired with the initial cluster $\{x_1,x_2,...,x_r\}$.
For the purposes of considering a consistent rank with enough interesting structure but not too large so as to require excessive computational resources all these algebras are of rank $r=4$.
Across the 3 types, 7 specific algebras were selected. 
These included 3 finite type algebras, generated from orientations of the $A_4$, $D_4$, and $F_4$ Dynkin diagrams (denoted A4, D4, F4 respectively). Note the $F_4$ algebra does not lie in the skew-symmetric classification, but is a skew-symmetrizable finite type. 
Next, 2 orientations of the affine $\tilde{A}_3$ type were used as finite-mutation type algebras which are specifically \textit{not} finite type; such that one had 1 anticlockwise 3-path and 1 clockwise arrow (denoted A13), whilst the other had a 2-path in each direction (denoted A22). 
Lastly, 2 infinite type algebras (denoted $\mathcal{I}$1 and $\mathcal{I}$2 respectively) were generated from the exchange matrices 

\begin{equation}
    EM_{\mathcal{I}1} = \begin{pmatrix}0 & 2 & 0 & 0\\ -2 & 0 & 1 & 0\\ 0 & -1 & 0 & 1\\ 0 & 0 & -1 & 0\end{pmatrix}\;,\qquad EM_{\mathcal{I}2} = \begin{pmatrix}0 & 2 & 0 & -2\\ -2 & 0 & 2 & 0\\ 0 & -2 & 0 & 1\\ 2 & 0 & -1 & 0\end{pmatrix}\;.
\end{equation}
The quivers for each of these initial seed exchange matrices are shown in Figure \ref{fig:quivers}.

Beyond the choice of Dynkin or affine Dynkin type, an orientation to each quiver must be prescribed.
The \texttt{sage} package initiates the finite type quivers with bipartite orientations, such that each node is either a source or a sink as shown.
Whilst the ambiguity to select an orientation may appear to lose generality, where the quiver's underlying graph is a tree (as for our finite types) any orientation is mutation equivalent to any other orientation \cite{CA_2}. 
Therefore for finite types, a choice of orientation is effectively a choice of initial point to expand around in the same exchange graph. 
Note this does not apply for the other mutation types, hence choosing any 2 different orientations produces 2 different cluster algebras and 2 different exchange graphs (as exemplified by the 2 chosen orientations of $\tilde{A}_3$ giving the 2 different A13 \& A22 algebras).
These cluster algebras were chosen such that under mutation similar Laurent polynomial styles (in particular monomial coefficients) were occurring so it was non-trivial to differentiate the seed representations by eye ahead of ML.

\begin{figure}[H]
	\centering
	\begin{subfigure}{0.45\textwidth}
    	\centering
    	\includegraphics[width=0.6\textwidth]{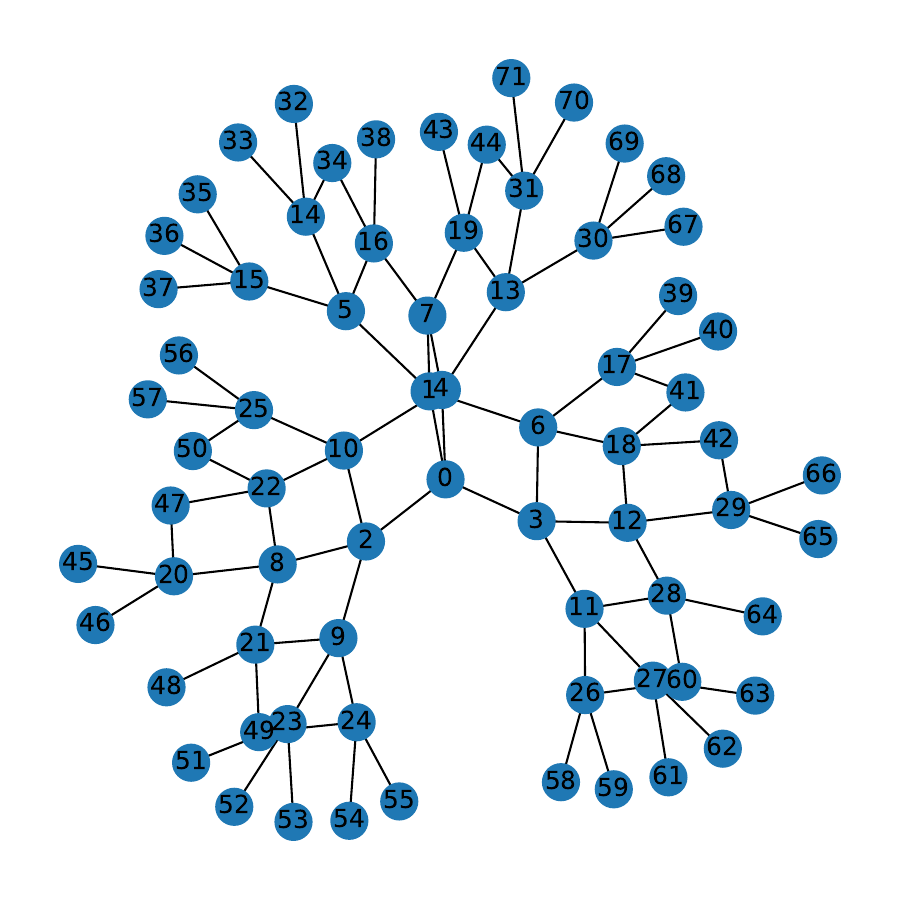}
    	\caption{A4 (finite)}\label{fig:A4EG}
	\end{subfigure} 
    \begin{subfigure}{0.45\textwidth}
    	\centering
    	\includegraphics[width=0.6\textwidth]{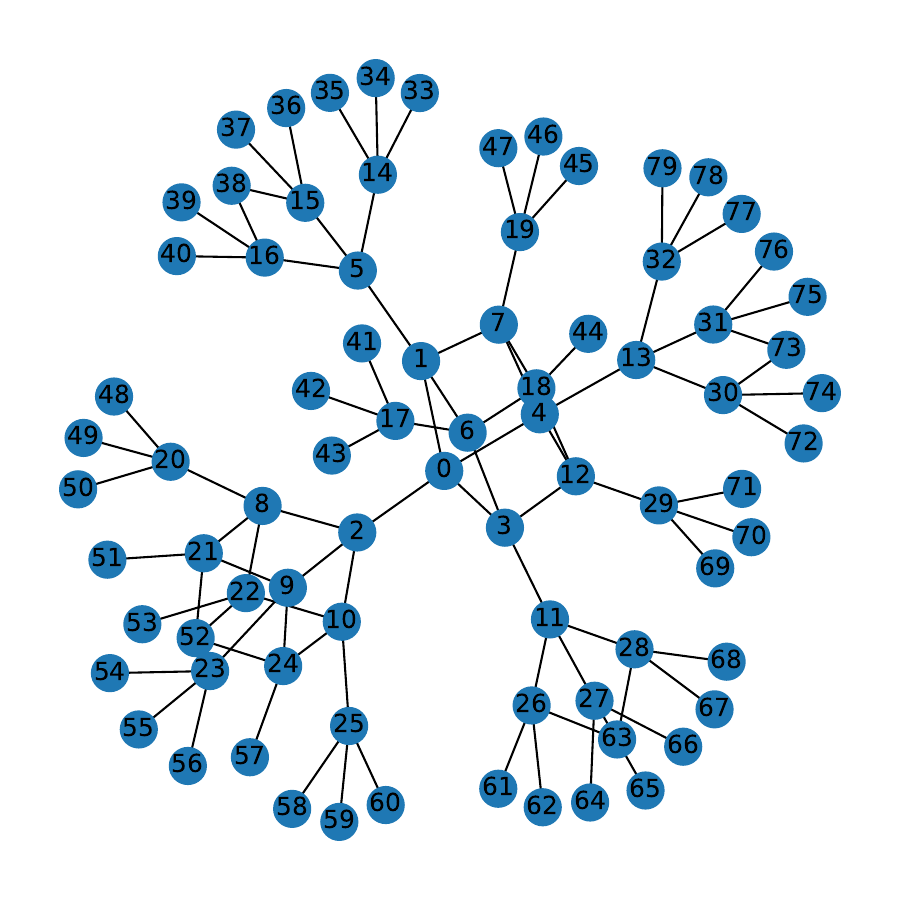}
    	\caption{D4 (finite)}\label{fig:D4EG}
    \end{subfigure} \\
    \begin{subfigure}{0.45\textwidth}
    	\centering
    	\includegraphics[width=0.6\textwidth]{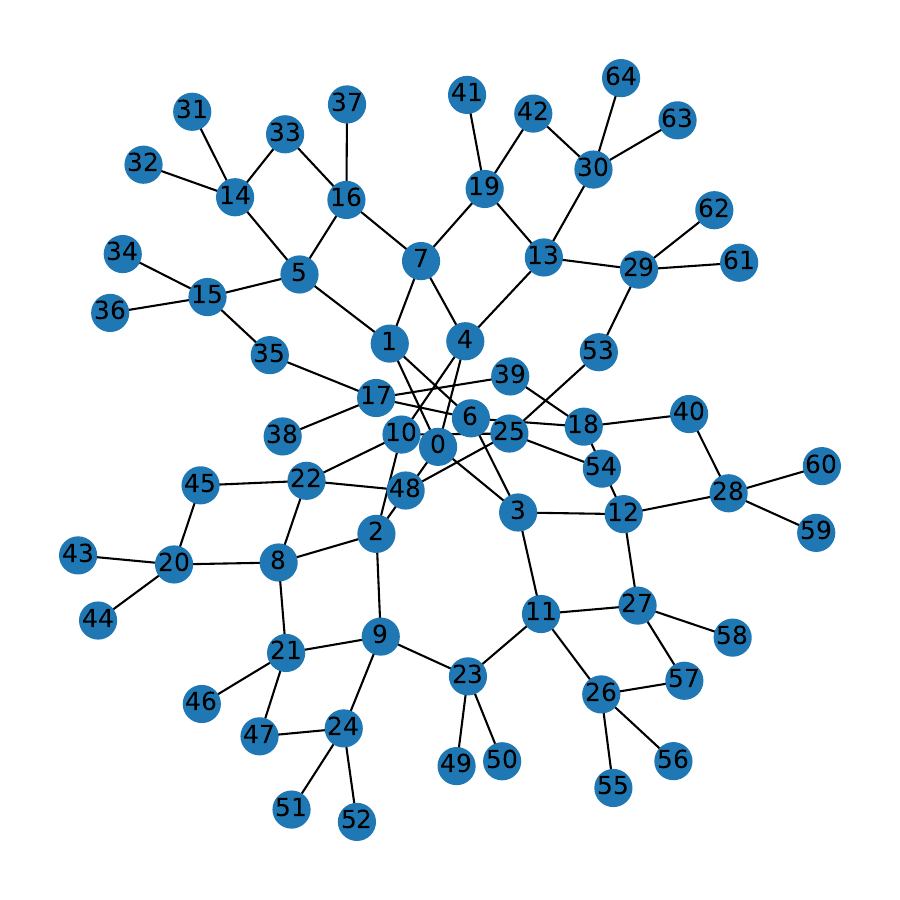}
    	\caption{F4 (finite)}\label{fig:F4EG}
    \end{subfigure} \\
	\begin{subfigure}{0.45\textwidth}
    	\centering
    	\includegraphics[width=0.6\textwidth]{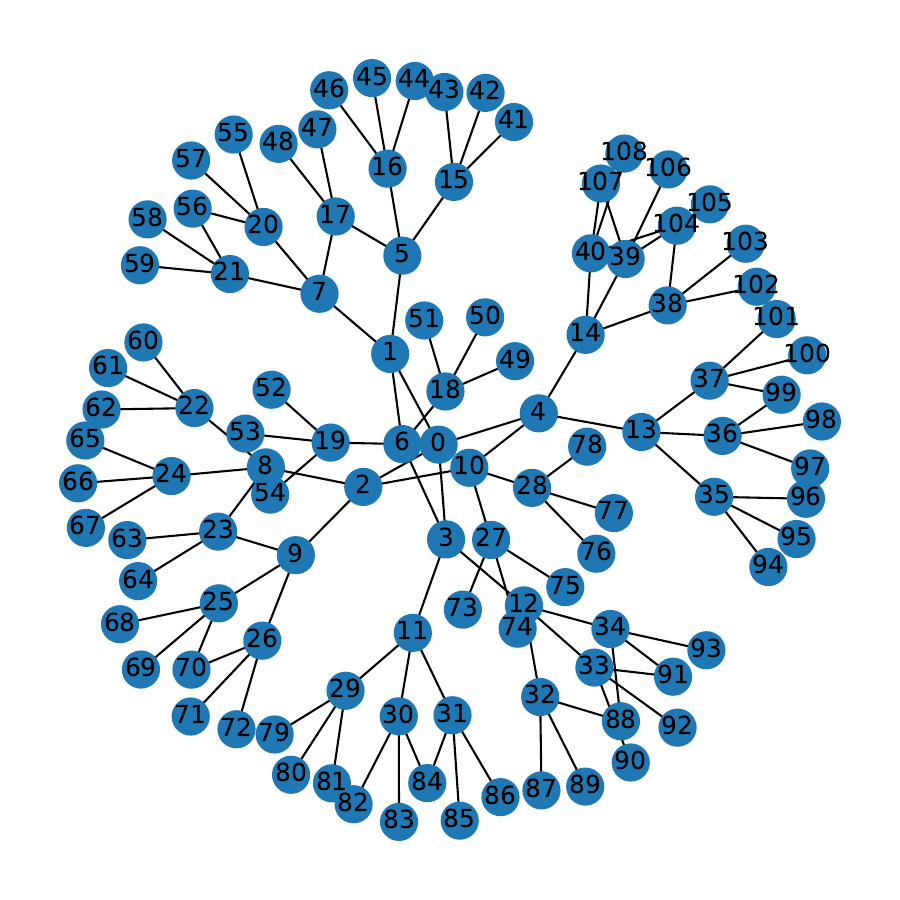}
    	\caption{A13 (finite-mutation)}\label{fig:A13EG}
	\end{subfigure} 
    \begin{subfigure}{0.45\textwidth}
    	\centering
    	\includegraphics[width=0.6\textwidth]{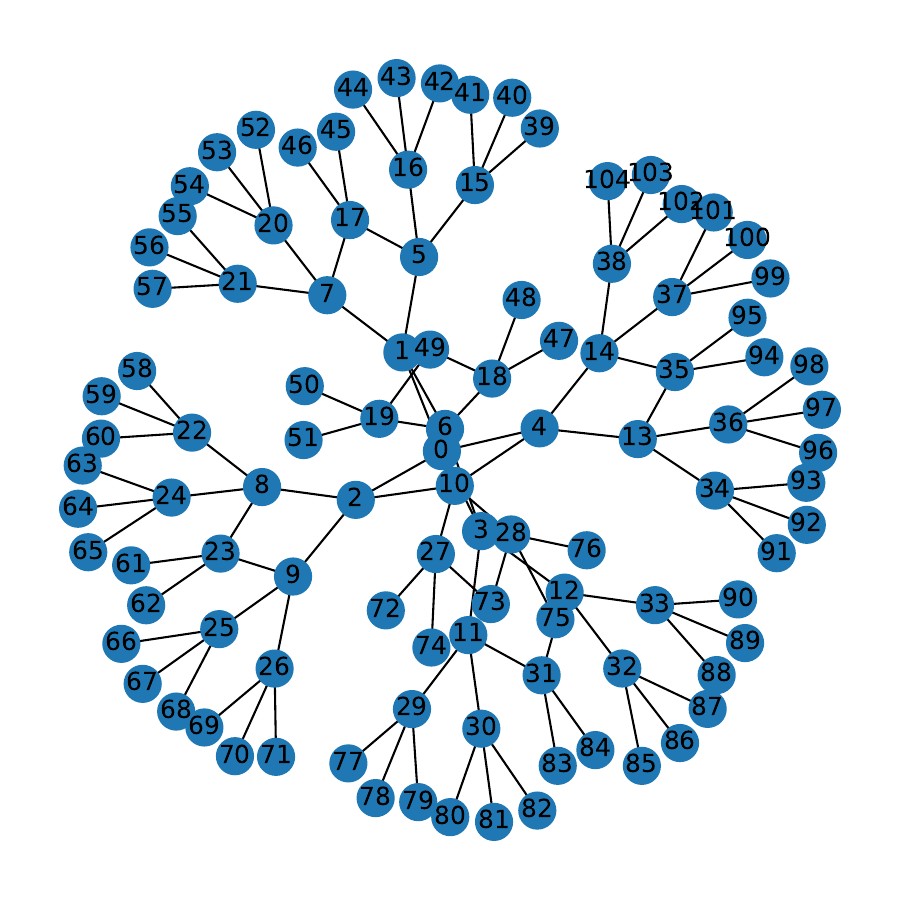}
    	\caption{A22 (finite-mutation)}\label{fig:A22EG}
    \end{subfigure} \\
	\begin{subfigure}{0.45\textwidth}
    	\centering
    	\includegraphics[width=0.6\textwidth]{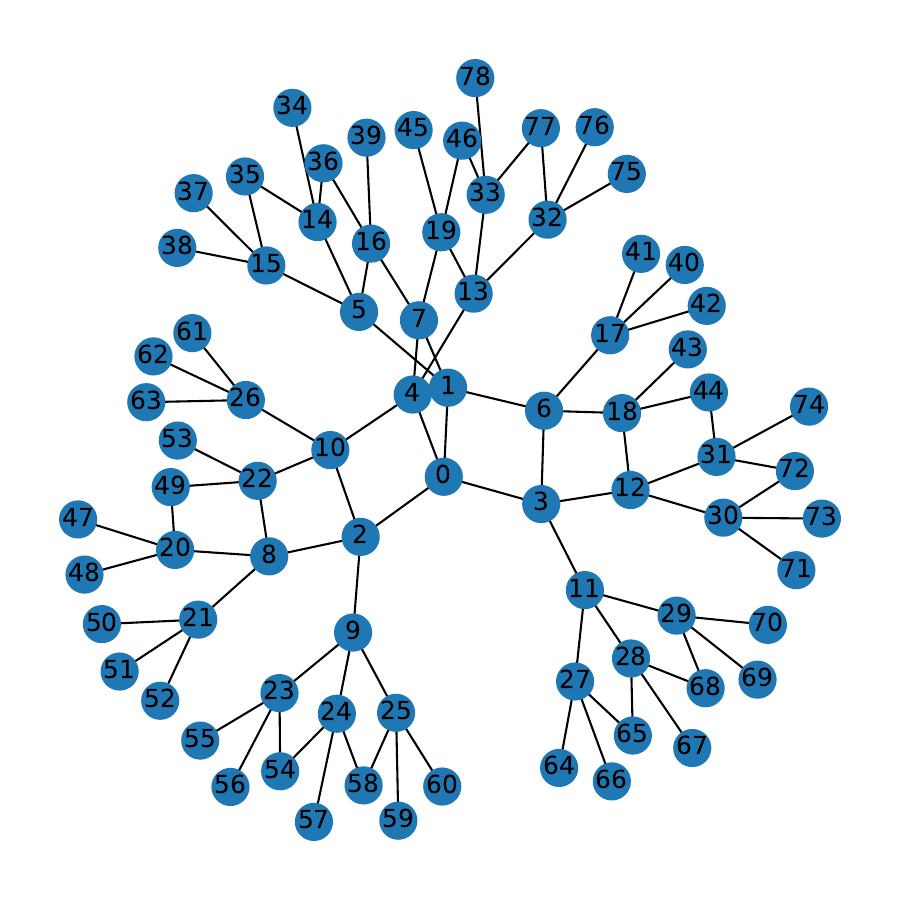}
    	\caption{$\mathcal{I}1$ (infinite)}\label{fig:inf1EG}
	\end{subfigure} 
    \begin{subfigure}{0.45\textwidth}
    	\centering
    	\includegraphics[width=0.6\textwidth]{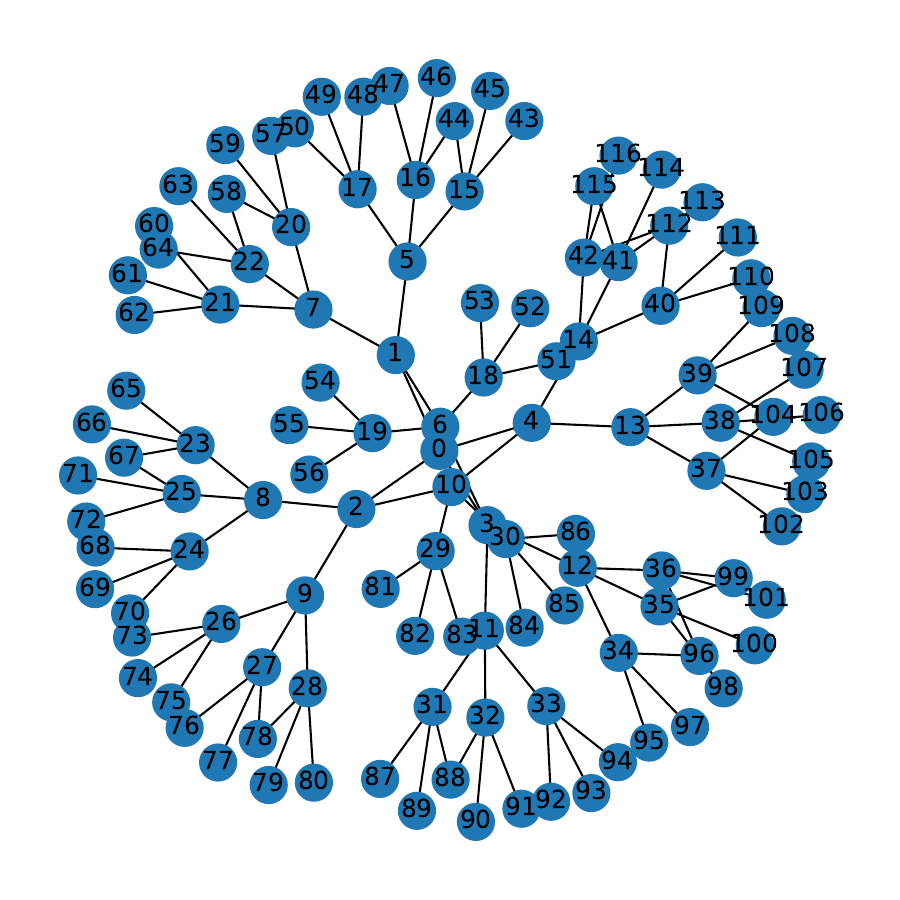}
    	\caption{$\mathcal{I}2$ (infinite)}\label{fig:inf2EG}
    \end{subfigure} \\
\caption{The seed exchange graphs generated to depth 4 for each of the considered cluster algebras. Types are labelled, where finite-mutation are specifically not finite type so are infinite for these seed exchange graphs but finite for the respective quiver exchange graphs (not shown). Vertices are labelled in the order they are generated starting from the initial seed `0'.}\label{fig:EGs}
\end{figure}

Due to the significant growth of complexity in the infinite type Laurent polynomials with depth, especially $\mathcal{I}2$ where depth 5 could not be computed in reasonable time, the core focus of the exchange graph analysis and subsequent ML was chosen to be up to and including depth 4.
In building the seed exchange graphs, clusters were not considered equivalent if their variables were the same but in a different order. As we see in later analysis identifying by this permutation equivalence loses some elegant symmetric structure in the exchange graphs. We account for the possibility of this causing some degeneracy in cycles on a case by case basis and mention where that happens explicitly.
Also, when defining an algebra only the cluster variables are important and hence taking the union of all clusters still produces the same generating set with or without application of this equivalence.
Exchange graphs up to depth 4 are given for each of the considered cluster algebras in Figure \ref{fig:EGs}.

\begin{figure}[!tb]
    \centering
    \includegraphics[width=0.6\textwidth]{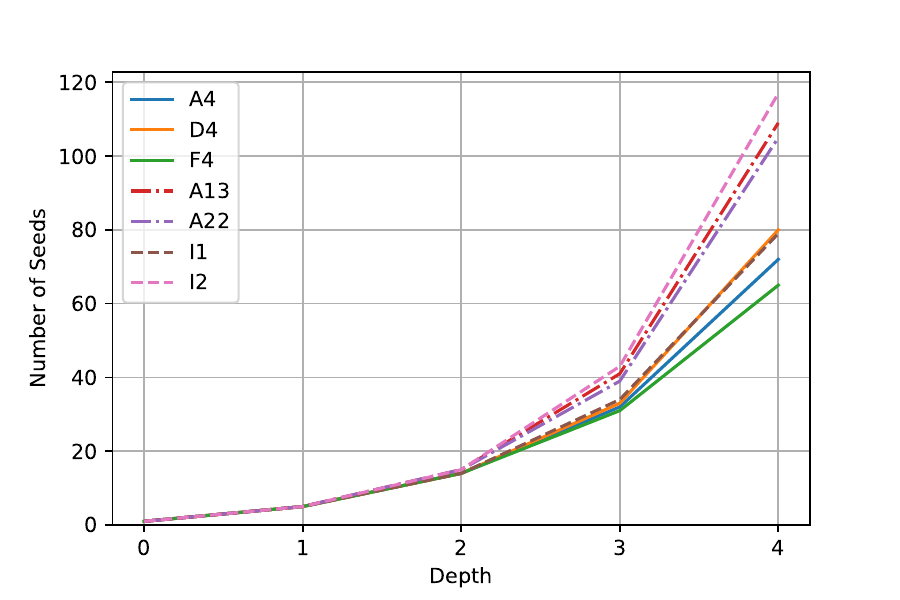}
    \caption{The number of seeds in the seed exchange graphs as depth varies for each of the considered cluster algebras, labelled by their respective initial seeds. Each type is depicted with a different linestyle.}
    \label{fig:EGclustswithdepth}
\end{figure}

Prior to any thorough data analysis of these exchange graphs, the types can begin to be distinguished 
by the number of seeds at each depth, as shown in Figure \ref{fig:EGclustswithdepth}, 
although it can be seen already that this information is not sufficient for classification: as more seeds are generated the finite type algebras are more likely to be reproducing previous seeds.
We believe that there is a similar behaviour where finite-mutation types (that are not finite) are also more likely to reproduce seeds than infinite types, as they are on the boundary between finite and infinite. 

To give an example for typical seed information we arbitrarily select seed 30 from the A4 seed exchange graph in Figure \ref{fig:A4EG}, whose exchange matrix and cluster take the form:
\begin{equation}
    EM_{A4:30} = \begin{pmatrix} 0 & -1 & 1 & 0 \\ 1 & 0 & -1 & 1 \\ -1 & 1 & 0 & 0 \\ 0 & -1 & 0 & 0 \end{pmatrix}\;,
\end{equation}
\begin{equation}
    \bigg\{x_1, \frac{(x_1x_3^2+x_1x_3+x_2x_4+x_3+1)}{x_2x_3x_4}, \frac{(x_2x_4+x_3+1)}{x_3x_4}, \frac{(x_3+1)}{x_4} \bigg\}\;.
\end{equation}

\subsection{Network Analysis}
The intricate tree structure of the exchange graphs exemplified by Fig.~\ref{fig:EGs} is suggestive for a treatment from network science.
For each of the cluster algebras considered the seed exchange graphs generated to depth 4 were thus compared with a variety of network analysis techniques.
The results of the analysis are provided in Table \ref{tab:EGanalysis}, and use a range of assessment techniques across the core themes of network analysis, including: clustering analysis, shortest path analysis, centrality analysis, and cycle basis analysis. 

\begin{table}[!tb]
\centering
\footnotesize
\begin{tabular}{|c|cccccc|}
\hline
\multirow{2}{*}{\begin{tabular}[c]{@{}c@{}}\\ Cluster\\ Algebra\end{tabular}} & \multicolumn{6}{c|}{Seed Exchange Graph Analysis (depth 4)} \\ \cline{2-7} 
& \multicolumn{1}{c|}{\begin{tabular}[c]{@{}c@{}}Number\\ of Vertices\end{tabular}} & \multicolumn{1}{c|}{Density} & \multicolumn{1}{c|}{\begin{tabular}[c]{@{}c@{}}Clustering\\ (tri, squ)\end{tabular}} & \multicolumn{1}{c|}{\begin{tabular}[c]{@{}c@{}}Wiener Index\\ (full, norm)\end{tabular}} & \multicolumn{1}{c|}{\begin{tabular}[c]{@{}c@{}}Centrality\\ (centre, diff)\end{tabular}} & \begin{tabular}[c]{@{}c@{}}Min cycle basis\\ ([length, freq])\end{tabular} \\ \hline
A4 & \multicolumn{1}{c|}{72} & \multicolumn{1}{c|}{0.034} & \multicolumn{1}{c|}{(0, 0.058)} & \multicolumn{1}{c|}{(13968, 5.46)} & \multicolumn{1}{c|}{(0, 0.029)} & [4,17] \\ \hline
D4 & \multicolumn{1}{c|}{80} & \multicolumn{1}{c|}{0.029} & \multicolumn{1}{c|}{(0, 0.037)} & \multicolumn{1}{c|}{(17941, 5.68)} & \multicolumn{1}{c|}{(0, 0.037)} & [4,13] \\ \hline
F4 & \multicolumn{1}{c|}{65}  & \multicolumn{1}{c|}{0.040} & \multicolumn{1}{c|}{(0, 0.064)} & \multicolumn{1}{c|}{(10700, 5.14)} & \multicolumn{1}{c|}{(0, 0.031)} & [4,17], [6,3] \\ \hline
A13 & \multicolumn{1}{c|}{109} & \multicolumn{1}{c|}{0.020} & \multicolumn{1}{c|}{(0, 0.034)} & \multicolumn{1}{c|}{(35284, 5.99)} & \multicolumn{1}{c|}{(0, 0.054)} & [4,12] \\ \hline
A22 & \multicolumn{1}{c|}{105} & \multicolumn{1}{c|}{0.021} & \multicolumn{1}{c|}{(0, 0.016)} & \multicolumn{1}{c|}{(32664, 5.98)} & \multicolumn{1}{c|}{(0, 0.061)} & [4,8] \\ \hline
$\mathcal{I}1$ & \multicolumn{1}{c|}{79} & \multicolumn{1}{c|}{0.031} & \multicolumn{1}{c|}{(0, 0.065)} & \multicolumn{1}{c|}{(17174, 5.57)} & \multicolumn{1}{c|}{(0, 0.015)} & [4,18] \\ \hline
$\mathcal{I}2$ & \multicolumn{1}{c|}{117} & \multicolumn{1}{c|}{0.019} & \multicolumn{1}{c|}{(0, 0.037)} & \multicolumn{1}{c|}{(41160, 6.07)} & \multicolumn{1}{c|}{(0, 0.063)} & [4,12] \\ \hline
\end{tabular}
\caption{Network analysis of the seed exchange graphs (EGs) generated to depth 4 for the 7 cluster algebras considered. The first 3 algebras are finite type, the latter 2 are infinite type, and the remaining 2 are finite-mutation type but not finite type (hence having infinite seed EGs). The analysis lists: the number of vertices in the EG; the density of the EG; the triangle and square average clustering coefficients; the Wiener index (both full form and normalised form); the eigenvector centrality analysis listing the central vertex and then the size of the smallest difference in centrality from the initial seed ``0'' to the clusters at depth 1; and finally the information on the minimum cycle basis showing the length of the basis cycles and the frequency of those lengths in the basis.  }\label{tab:EGanalysis}
\end{table}

For each exchange graph (sometimes denoted EG) the number of distinct seeds (i.e. vertices in the EG) up to depth 4 is given.
As expected the infinite cluster algebras (including A13 and A22 which are distinctly not finite) have more total vertices, except for $\mathcal{I}1$ which is surprisingly low.
This is due to the finite type algebras usually reproducing more previously generated seeds as mutation continues to higher depths.

Further vertex analysis usually considers the degree distribution, but due to the construction process for seed exchange graphs all vertices will have degree 4, except those truncated from further mutation by the depth limit.
Therefore this would provide little insightful analysis for this graph style.

The exchange graph \textit{density} then considers the number of edges, as opposed to the number of vertices. 
Here, the total number of edges in the graph is divided by the total number of possible edges (i.e. the number of edges in a complete graph with the same number of vertices).
Due to there being 4 edges per vertex in the depths up to 3, as each cluster has 4 variables to mutate, these density scores correlate with the number of vertices.
The finite types have higher densities as they are more tightly-knit graphs and are further from the more tree-like structure of the infinite types, as can be seen in the graphs of Figure \ref{fig:EGs}.
In addition, the finite-mutation types with infinite exchange graphs have a slightly higher density than most infinite types (better represented by $\mathcal{I}2$), supporting that they are on this border of the seed exchange graphs closing up from infinite to finite. 
Although $\mathcal{I}1$ has a higher density, closer to the finite types at this depth, as can be seen in its seed exchange graph, it does show the typical tree-like substructure of infinite types, with an unusual star-pattern of lines of 4-cycles forming a net-like structure between the tree subgraphs coming off them.

\paragraph{Clustering Analysis}
\textit{Clustering coefficients} give information about how vertices cluster within graphs.
The two styles of coefficient considered here are both global in nature, and consider the number of triangles (3-cycles) or squares (4-cycles) that exist in the graph relative to the total number of possible triangles or squares that could exist.

For all the algebras there were no triangles in the seed exchange graphs. 
This is expected since it would require seeds to either jump mutation depths or connect the same depth.
Both these scenarios are not possible since any 2-path of connected vertices either spans 3 depths so closing this 2-path into a 3-cycle would require the mutation that this closing edge represents to jump a depth (i.e. mutating two variables simultaneously), or the 2-path has the initial seed as its centre vertex and hence the vertices which need to connect would both be depth 1 and unable to mutate to each other, as they are both a different mutated-variable from the initial seed, hence again requiring a double mutation.

This idea generalises to all odd size cycles. 
Since each mutation changes the depth of seed under consideration, for any sequence of mutations to close into a cycle the number of mutations increasing depth must equal the number reducing the depth. 
Therefore all cycles must be of even length. 

However all the algebras have a non-zero square clustering coefficient.
Squares, or 4-cycles, in a seed exchange graph indicate the scenario where two mutation actions commute; i.e. a seed can have two variables mutated in either order to produce the same seed.
The frequency of this commutative action interestingly does not appear to correlate with the algebra type, exemplified by D4 having the same coefficient as $\mathcal{I}2$ and A13 having over double the coefficient value of A22.

\paragraph{Shortest Path Analysis}
The shortest path analysis carried out comes in the form of the \textit{Wiener index}.
This index computes the sum of the shortest paths between all pairs of vertices (the full form), which we also normalise by the number of pairs of vertices, $nC2$, to give an average shortest path between vertices.
The normalised form is more useful for comparison and indicates that infinite types (excluding $\mathcal{I}1$) have vertices further separated on average.
This measure provides an indication of the frequency of cycles, and also the placement of them. 
If they are more spread over the graph (as for F4) they are more useful in shortening the shortest paths between outer vertices.

\paragraph{Centrality Analysis}
Centrality of a network determines the natural centre vertex. 
Due to the generation process being from an initial seed to some depth, this would make the initial seed (labelled 0) the logical choice. 
However, experimentation with some lower rank algebras showed this was not always the case so it is good to confirm for these  algebras considered.
The centrality measure used was \textit{eigenvector centrality}, which performs eigendecomposition of the graph's adjacency matrix, then taking the eigenvector with the largest eigenvalue, whose entries indicate the relative importance of each vertex in the network.
In addition to the computed centre vertex (highest respective entry), the difference in centrality between the centre initial seed and the most central seed at depth 1 is given.
Where this difference is large, as for most of the infinite types, the initial seed is a more obvious centre.
However, the difference is remarkably small for $\mathcal{I}1$, which we believe to be attributed to the symmetry of 4-cycles surrounding the centre cycle (with a narrow diamond shape) with vertices labelled from [0,4,7,1].

\paragraph{Cycle Basis Analysis}
Where the exchange graphs diverge from the limiting case of an infinite tree, cycles are introduced to the graph as seeds are reproduced through mutation. 
As more cycles are added the graph begins to close, becoming `more finite' in the process with less new seeds produced per mutation.
It is hence this cycle structure we believe to be key to what dictates the cluster algebra behaviour, correlating loosely with the algebra type.
Not only the frequency of cycles, but their distribution of sizes is an area of particular interest.

To analyse the cycle structure of each seed exchange graph, the \textit{minimum cycle basis} for each was computed.
The minimum cycle basis puts emphasis on selecting basis cycles with the lowest lengths, which can then be summed by symmetric difference in the cycle vector space over the finite field of 2 elements, $\mathbb{Z}/2\mathbb{Z}$, to produce any cycle in the graph.
Since the symmetric difference of any 2 even cycles is also even ($s$-cycle + $t$-cycle with $u$ overlapping edges = $(s+t-2u)$-cycle), this corroborates the idea that all cycles being even is self-consistent.
A corollary of this is that all cycles in any basis will be even also.
In Table \ref{tab:EGanalysis}, the structure of the minimum cycle bases for each algebra is given.
Since each cycle space can have many legitimate bases, the actual basis content was not of focus; instead the number of cycles required (`freq') and of what sizes (`length') provide the interesting information for analysis.

It turns out for most of these algebras, that the cycle spaces at this depth can be constructed from exclusively 4-cycles (this does not hold in general at higher depths).
4-cycles are the natural lower bound of cycle sizes in seed exchange graphs.
This is because any mutation sequence (longer than the trivial involution 2-cycle for all edges) has to mutate on any vertex at least twice to change it back, hence 4-cycles being the minimum where two mutations on different vertices commute.
However, F4 requires 3 6-cycles to define the cycle space up to depth 4, implying that its exchange graph structure is distinct from the others with a more subtle design.
This may be related to the skew-symmetrizable vs skew-symmetric style of the exchange matrices, but more investigation into this would be required for a range of skew-symmetrizable matrices before concluding any particular pattern.

Interestingly the finite-mutation types A13 \& A22 have the smallest frequency of 4-cycles in their basis:
expectedly less than the finite type, since their graphs are compact, however also less than the infinite types.
This may be due to a more systematic occurrence of cycles, which are further apart and hence fewer are visible at this depth; whereas the infinite graphs have a skewed distribution of cycles where there is some iterative commutativity via 4-cycles.

\begin{figure}[!tb]
	\centering
	\begin{subfigure}{0.45\textwidth}
    	\centering
    	\includegraphics[width=0.9\textwidth]{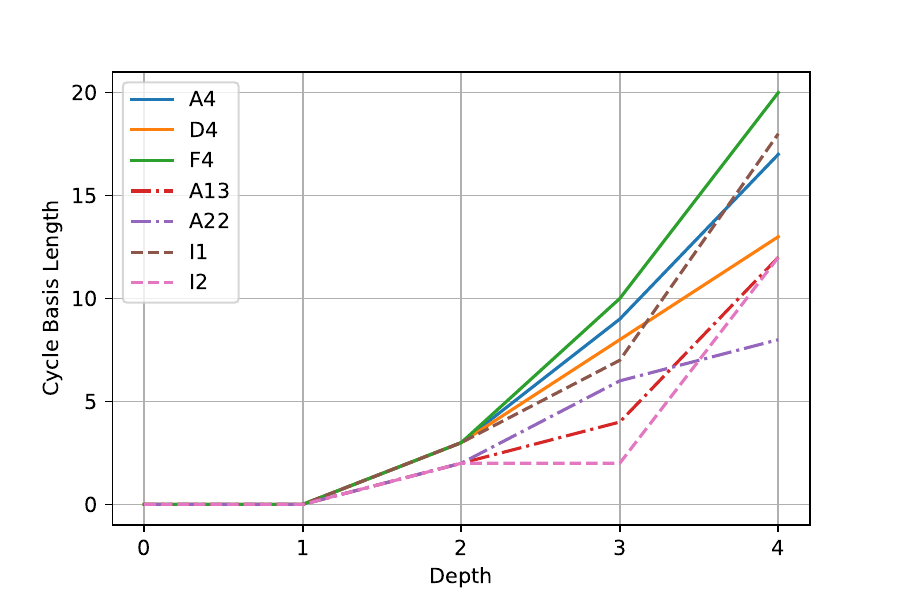}
    	\caption{Cycle basis size}\label{fig:cyclebasis}
	\end{subfigure} 
    \begin{subfigure}{0.45\textwidth}
    	\centering
    	\includegraphics[width=0.9\textwidth]{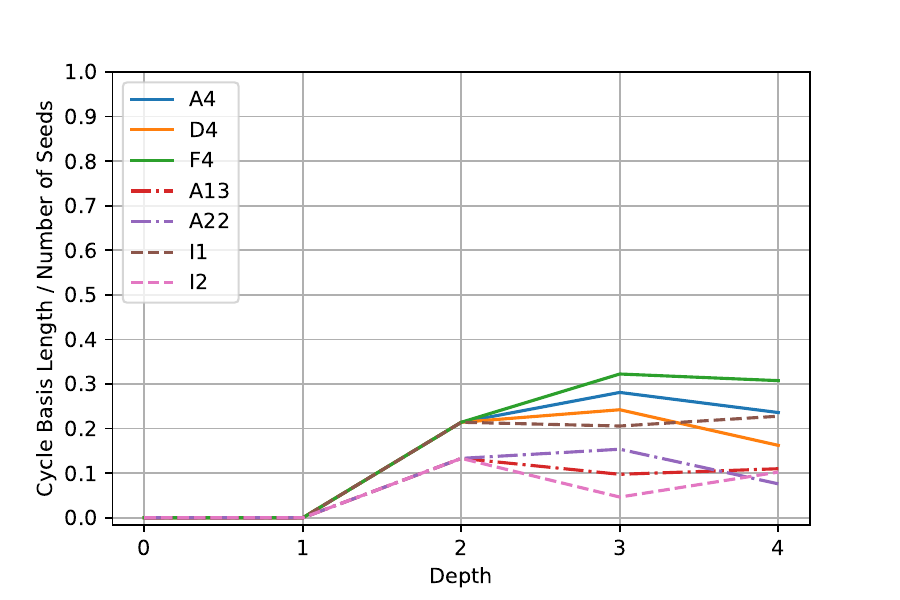}
    	\caption{Relative cycle basis size}\label{fig:relativecyclebasis}
    \end{subfigure} \\
    \begin{subfigure}{0.45\textwidth}
    	\centering
    	\includegraphics[width=0.9\textwidth]{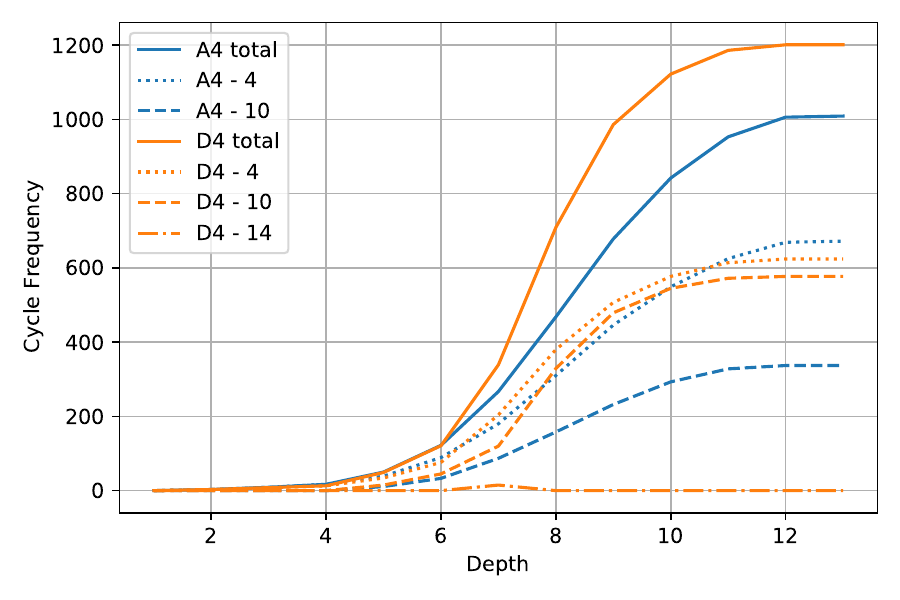}
    	\caption{A4 \& D4 cycle bases sizes}\label{fig:A4D4cyclebasis}
    \end{subfigure} \\
\caption{The size of the minimum cycle basis for each of the considered cluster algebras' seed exchange graphs as depth varies is plotted directly in (a), and the cycle basis length relative to the number of seeds is plotted in (b). Moreover, (c) focuses on the finite type A4 and D4 algebras with seed exchange graphs generated to their maximum depths; beyond depth 4 10-cycles are introduced so each cycle length frequency is denoted separately, as well as the total frequencies.}\label{fig:CycleBasisLengthwithDepth}
\end{figure}

How the lengths of the cycle bases vary with depth (up to the computable depth 4 for infinite type, and beyond for the finite types) is given in Figure \ref{fig:CycleBasisLengthwithDepth}.
Since all basis cycles were of length 4 except for 3 6-cycles in the F4 case, the total number of cycles in the basis is plotted instead of differentiating these cycle lengths.
It shows that there is some noise of basis length growth at these lower depths for the algebras, and thus particularly for the infinite types it would be useful to find more efficient ways of computing the seed exchange graphs to probe higher depths.
The basis length with respect to the number of seeds is also shown; since behaviour is similar to without this normalisation it indicates that there is some correlation described by more seeds at a certain depth indicating the seeds are less connected leading to a smaller relative basis size.

Interestingly, the A4 and D4 finite type algebras have similar behaviour up to depth 5. 
This is the depth at which 10-cycles are permitted, and beyond this depth the behaviour diverges with D4 having consistently more total cycles due to the considerably larger number of 10-cycles.
This indicates that these algebras have fundamentally different construction; A4's more significant structural reliance on 4-cycles may be due to more mutation operations commuting which in turn is likely due to the lower connectivity of the original A4 quiver, translating to subsequent quivers also.
We note there is an anomalous 15 14-cycles occurring for D4 at depth 7: an artefact of the depth truncation that immediately splits into 4- and 10-cycles at the next depth.

Additionally, dual to the cycle vector space is the cut vector space; and although an interesting alternative method of network analysis, we leave its application to exchange graphs to future work.

\subsection{Generalised Associahedra \& Seed Equivalence}
The cluster variables for the infinite type algebras quickly become highly complex with mutation -- this prevents generation of the seed exchange graphs to large depths.
However the finite type algebras do not suffer as severely from this behaviour, and due to their finite-ness their entire exchange graphs (`generalised associahedra') of all seeds can be completely generated.

\begin{figure}[!tb]
	\centering
	\begin{subfigure}{0.45\textwidth}
    	\centering
    	\includegraphics[width=0.9\textwidth]{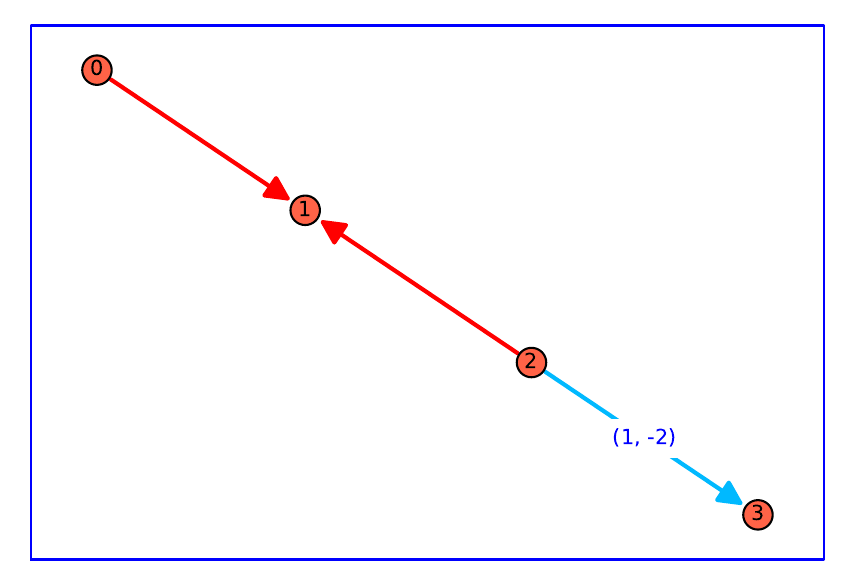}
    	\caption{B4}\label{fig:B4quiv}
	\end{subfigure} 
    \begin{subfigure}{0.45\textwidth}
    	\centering
    	\includegraphics[width=0.65\textwidth]{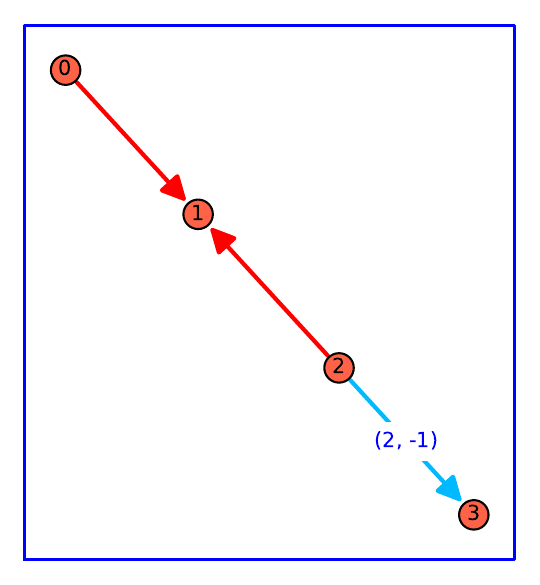}
    	\caption{C4}\label{fig:C4quiv}
    \end{subfigure} \\
\caption{Quivers defining the exchange matrices for the remaining rank 4 finite type cluster algebra initial seeds. Both (a) B4 and (b) C4 are skew-symmetrizable, shown by the non-simply laced, double weighted, edges of opposite weighting.}\label{fig:B4C4Quivers}
\end{figure}

Therefore for the rank 4 finite type cases we generate their generalised associahedra, and perform similar network analysis on them all.
Beyond the A4, D4, F4 cases considered previously, we also introduce the B4 and C4 cluster algebras.
These are also finite type, and as they arise from non-simply laced Dynkin orientations, they are not skew-symmetric but skew-symmetrizable. 
This Dynkin terminology means simply-laced edges have weight one, whilst the non-simply laced edges have a double weighting given by the non-skew-symmetric exchange matrix components.
Their respective quivers are shown in Figure \ref{fig:B4C4Quivers}, and their network analysis along with A4, D4, \& F4 in Table \ref{tab:FinitefullEGanalysis}.

Interestingly, all the skew-symmetrizable quivers generated algebras with the same number of seeds and density at the same depth.
Further to this, B4 and C4 have identical analysis values for all measures, indicating the swap of quiver edge weighting separating them has no effect on the seed exchange graph.
F4 has some differences, with slightly fewer squares but a better connectivity (via Wiener index), supported by more 6-cycles and fewer 4-cycles in the basis.
The B4 and C4 algebras actually have the same seed exchange graph; however the Laurent polynomials of the variables take different values and hence the seeds at each vertex are very different \cite{felikson2010unfold,musiker2011compendium}.

Conversely, the A4 and D4 algebras have similar numbers of clusters and densities, despite D4 having far fewer 4-cycles and many more 10-cycles in its basis.
Since D4's fundamental structure relies more on 10-cycles, these provide quicker paths between components of the generalised associahedra leading to a smaller normalised Wiener index.

Gratifyingly, all these generalised associahedra have no discernible centre via eigenvector centrality, such that the dominant eigenvector has all its entries equal for each algebra.
This supports that any seed may be used to generate each algebra symmetrically, as these generalised associahedra form complicated polytopes built from 4- and 10-cycles; and in the skew-symmetrizable cases 6-cycles become a necessity too. 

\begin{table}[!tb]
\centering
\addtolength{\leftskip} {-2cm}
\addtolength{\rightskip}{-2cm}
\footnotesize
\begin{tabular}{|c|cccccc|}
\hline
\multirow{2}{*}{\begin{tabular}[c]{@{}c@{}}\\ Cluster\\ Algebra\end{tabular}} & \multicolumn{6}{c|}{Generalised Associahedra Analysis} \\ \cline{2-7} 
& \multicolumn{1}{c|}{\begin{tabular}[c]{@{}c@{}}Number\\ of Vertices\end{tabular}} & \multicolumn{1}{c|}{Density} & \multicolumn{1}{c|}{\begin{tabular}[c]{@{}c@{}}Clustering\\ (tri, squ)\end{tabular}} & \multicolumn{1}{c|}{\begin{tabular}[c]{@{}c@{}}Wiener Index\\ (full, norm)\end{tabular}} & \multicolumn{1}{c|}{\begin{tabular}[c]{@{}c@{}}Centrality\\ (centre vertex)\end{tabular}} & \begin{tabular}[c]{@{}c@{}}Min cycle basis\\ ([length, freq])\end{tabular} \\ \hline
A4 & \multicolumn{1}{c|}{1008 (13)} & \multicolumn{1}{c|}{0.0040} & \multicolumn{1}{c|}{(0, 0.080)} & \multicolumn{1}{c|}{(3881976, 7.65)} & \multicolumn{1}{c|}{*no centre*} & [4,672], [10,337] \\ \hline
B4 & \multicolumn{1}{c|}{420 (10)} & \multicolumn{1}{c|}{0.0095} & \multicolumn{1}{c|}{(0, 0.077)} & \multicolumn{1}{c|}{(542400, 6.16)} & \multicolumn{1}{c|}{*no centre*} & [4,270], [6,60], [10,91] \\ \hline
C4 & \multicolumn{1}{c|}{420 (10)}  & \multicolumn{1}{c|}{0.0095} & \multicolumn{1}{c|}{(0, 0.077)} & \multicolumn{1}{c|}{(542400, 6.16)} & \multicolumn{1}{c|}{*no centre*} & [4,270], [6,60], [10,91] \\ \hline
D4 & \multicolumn{1}{c|}{1200 (12)} & \multicolumn{1}{c|}{0.0033} & \multicolumn{1}{c|}{(0, 0.072)} & \multicolumn{1}{c|}{(5150592, 7.16)} & \multicolumn{1}{c|}{*no centre*} & [4,624], [10,577] \\ \hline
F4 & \multicolumn{1}{c|}{420 (10)} & \multicolumn{1}{c|}{0.0095} & \multicolumn{1}{c|}{(0, 0.072)} & \multicolumn{1}{c|}{(536816, 6.10)} & \multicolumn{1}{c|}{*no centre*} & [4,252], [6,111], [10,58] \\ \hline
\end{tabular}
\caption{Network analysis of the generalised associahedra for the finite type cluster algebras, labelled by their respective initial seeds. The analysis lists: the number of vertices in the EG (with depth to generate these in brackets); the density of the EG; the triangle and square average clustering coefficients; the Wiener index (both full form and normalised form); the eigenvector centrality analysis listing the central vertex (or lack of); and finally the information on the minimum cycle basis showing the length of the basis cycles and the frequency of those lengths in the basis.}\label{tab:FinitefullEGanalysis}
\end{table}

\paragraph{Seed Equivalence}
Since only the cluster variables are required to define a cluster algebra, once all variables have been generated the quivers and clusters themselves become superfluous information.
Therefore usually when one considers clusters one considers them equivalent up to permutation of the variables within the cluster, where the exchange matrix (and hence quiver) must be permuted in the same way. 
However, we find from some testing that different numbers (i.e. $\leq r!$) of permutations of clusters are produced within each algebra. 
Therefore by identifying under the full permutation group some combinatorial structure of the generation process is lost.

For example, for the generalised associahedra considered here, which the entire seed exchange graphs can be generated for, one expects the number of clusters (up to permutation equivalence), $N$, to follow the relation
\begin{equation}
    N=\prod_{i=1}^n \frac{e_i+h+1}{e_i+1}=\prod_{i=1}^n \frac{d_i+h}{d_i}\;,
\end{equation}
for $e_i$ the exponents and $d_i$ the degrees (of polynomial invariants) of the considered Dynkin type's root system (of which there are $n$), and $h$ the Coxeter number \cite{fomin2003systems,fomin2008root, Dechant2018Trinity}. 
For each of the rank 4 finite type algebras considered, the number of clusters up to permutation equivalence, $N$, is given in Table \ref{tab:permequivfactors}.
One may then naively expect all permutations of the seeds to occur in the exchange graphs generated when this permutation equivalence is not identified by.
However, although the A4 and D4 algebras do have all $4!=24$ permutations, the non-simply laced types do not.
The number of seeds without identification by permutation equivalence that we generate in the full generalised associahedra, denoted $N'$, are repeated in Table \ref{tab:permequivfactors} also for reference.

\begin{table}[!tb]
\centering
\begin{tabular}{|c|c|c|c|c|c|}
\hline
Cluster Algebra & A4   & B4  & C4  & D4   & F4  \\ \hline
$N$       & 42   & 70  & 70  & 50   & 105 \\ \hline
$N'$ & 1008 & 420 & 420 & 1200 & 420 \\ \hline
$N'/N$  & 24   & 6   & 6   & 24   & 4   \\ \hline
\end{tabular}
\caption{The rank 4 finite type cluster algebras considered whose generalised associahedra are generated in full. $N$ is the number of seeds up to identification by the permutation equivalence, $N'$ the number of seeds generated without the permutation equivalence, and the final row the factor between them.}\label{tab:permequivfactors}
\end{table}

As exemplified by the rank 4 cases, the non-simply laced cluster algebras (B4, C4, F4) do not generate all permutations of seeds from mutation about an initial seed.
In fact all permutations are allowed amongst the simply-laced components, where components are defined to be the sets of quiver nodes (and hence cluster positions) only connected by simply-laced edges, i.e. components are connected by non-simply laced edges.
Therefore any cluster $\{x_1,x_2,x_3,x_4\}$ will be mutation equivalent to $\{x_1,x_3,x_2,x_4\},\{x_2,x_1,x_3,x_4\},\{x_2,x_3,x_1,x_4\},\{x_3,x_1,x_2,x_4\},\{x_3,x_2,x_1,x_4\}$ for the B4 and C4 algebras; whilst mutation equivalent to $\{x_1,x_2,x_4,x_3\},\{x_2,x_1,x_3,x_4\},\\ \{x_2,x_1,x_4,x_3\}$ for the F4 algebra.
We emphasise here that this behaviour holds for all clusters of any variables (i.e. including those which are Laurent polynomials of the initial variables), hence occurring with permutation frequencies given by the relevant factor.

Beyond the rank 4 cases considered we find through experimentation up to rank 5 that the factor $N'/N$ is $r!$ for the $A_r$ and $D_r$ types, $(r-1)!$ for the $B_r$ and $C_r$ types (anticipating these to also hold $\forall \ r)$, and $1$ for $G_2$ (where all 8 clusters are different combinations of cluster variables). 
We therefore predict due to the lack of non-simply laced edges that the $E_r$ types have an $r!$ factor too (for $r\in\{6,7,8\}$), despite them being too large to computationally generate in full.

This factor for the non-simply laced types is related to their skew-symmetrizable exchange matrices. 
Since the skew-symmetrizer diagonal matrix in the skew-symmetrization process will have a non-unit factor associated to the non-simply laced edge, and since the skew-symmetrizer matrix is preserved under mutation \cite{nakanishi2022cluster}, any mutation involving a non-simply laced edge will cause the variables crossing it to pick up some non-trivial power that cannot be cancelled in this other component.
Therefore any cluster variable is restricted in its current form to only appear in its current component (separated by a non-simply laced edge to other components).
We also note that whilst mutation may introduce more non-simply laced edges, these will only ever connect nodes across the different components, and hence the above reasoning still applies for permutations of clusters within the components, which are dictated by the initial seed.

This subtle combinatorial structure of cluster mutation and the algebra generation process is lost under the identification by permutation equivalence.
Since our primary focus here is the exchange graphs built from these clusters, we think it best to consider the full structure for each algebra without the permutation equivalence applied, with the idea that once this full exchange graph is generated one can then still take all independent variables from all clusters (which already have a lot of overlap) to retrieve the algebra's generators.
This viewpoint also reveals other unanticipated structure in the exchange graphs as detailed in the following section.

\subsection{Quiver Exchange Graphs}
Whilst the previous subsection analysed the seed exchange graphs, this subsection shifts focus to the smaller quiver exchange graphs, with hope of uncovering the embedding structure into the seed exchange graph form (where the cluster information is included).

\begin{figure}[H]
	\centering
	\begin{subfigure}{0.45\textwidth}
    	\centering
    	\includegraphics[width=0.6\textwidth]{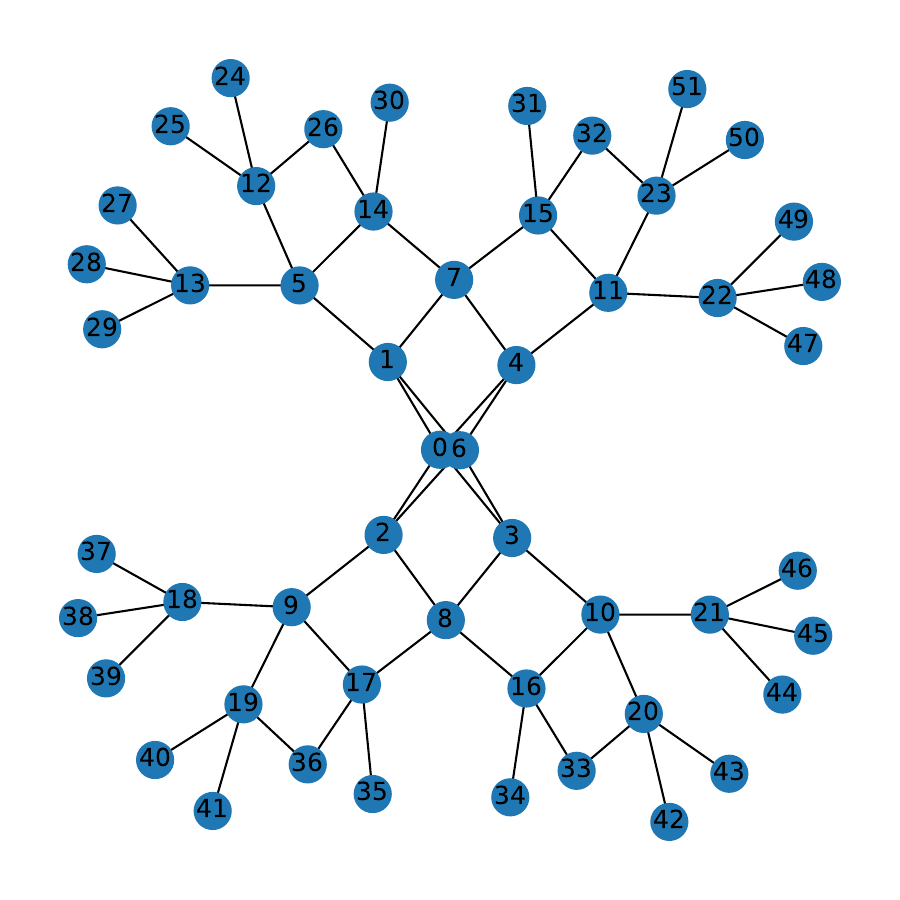}
    	\caption{A4 (finite)}\label{fig:A4EG}
	\end{subfigure} 
    \begin{subfigure}{0.45\textwidth}
    	\centering
    	\includegraphics[width=0.6\textwidth]{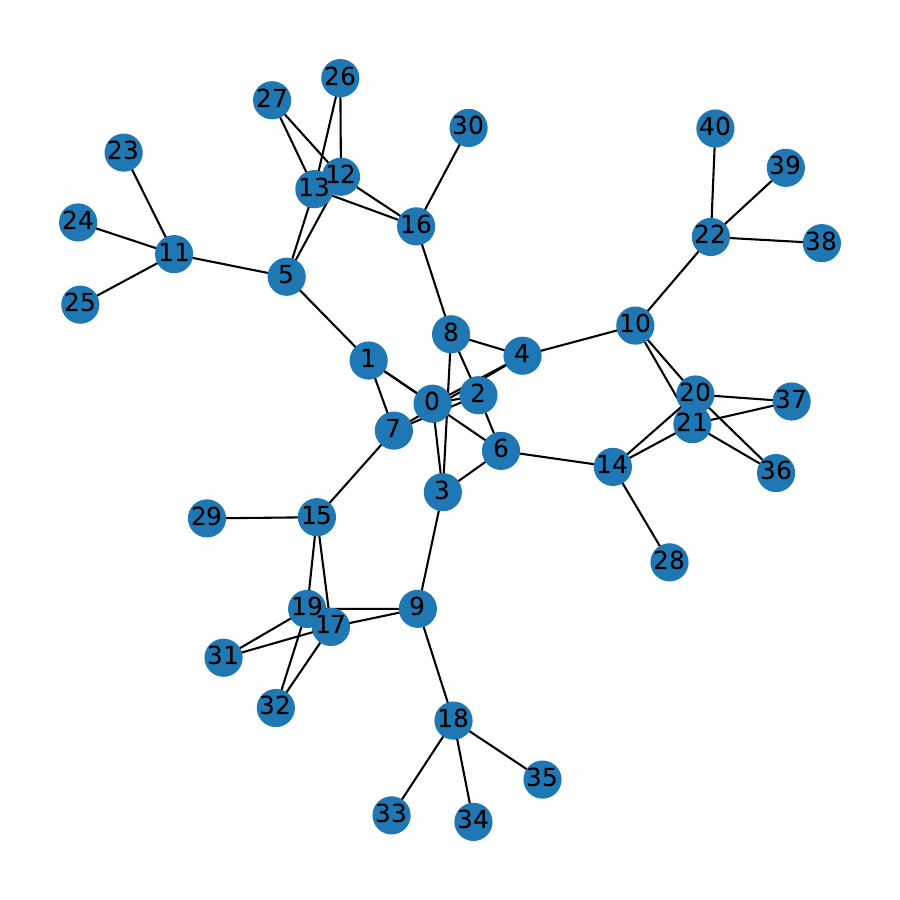}
    	\caption{D4 (finite)}\label{fig:D4EG}
    \end{subfigure} \\
    \begin{subfigure}{0.45\textwidth}
    	\centering
    	\includegraphics[width=0.6\textwidth]{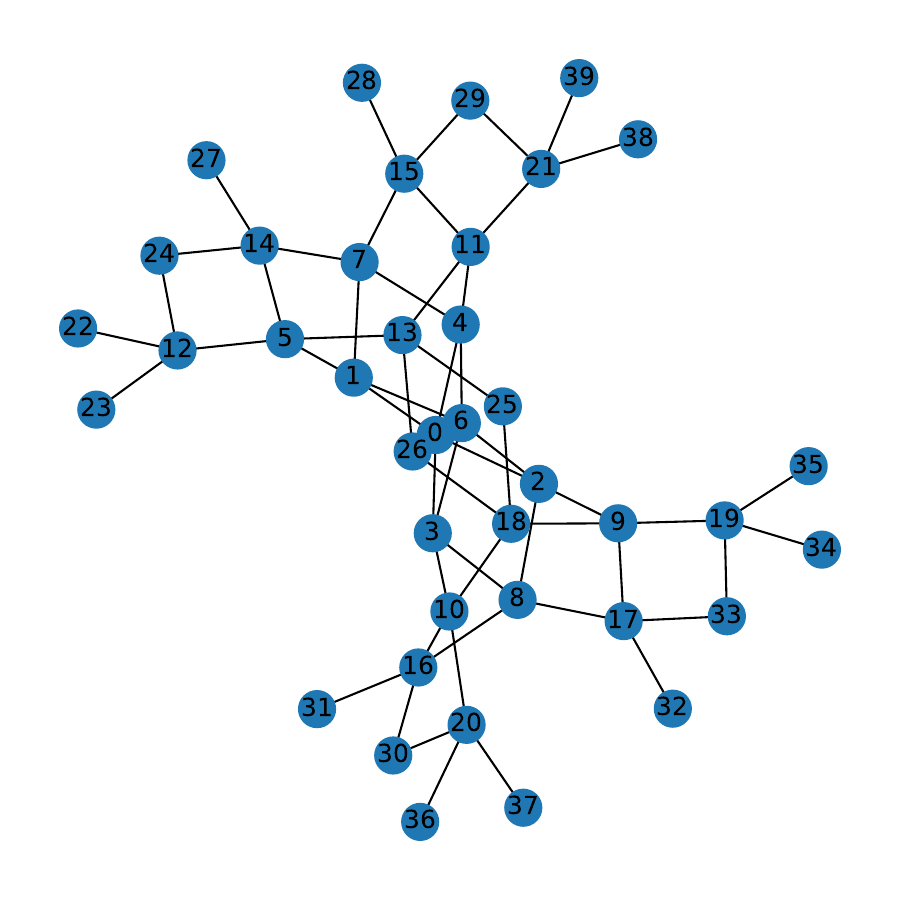}
    	\caption{F4 (finite)}\label{fig:F4EG}
    \end{subfigure} \\
	\begin{subfigure}{0.45\textwidth}
    	\centering
    	\includegraphics[width=0.6\textwidth]{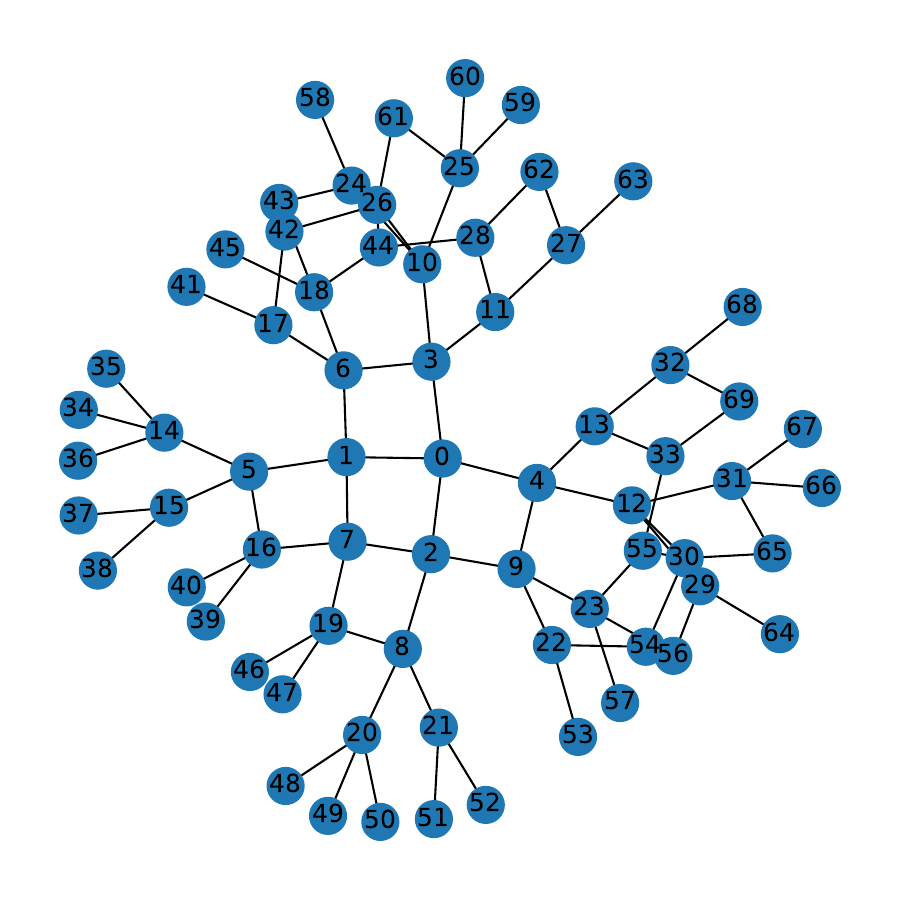}
    	\caption{A13 (finite-mutation)}\label{fig:A13EG}
	\end{subfigure} 
    \begin{subfigure}{0.45\textwidth}
    	\centering
    	\includegraphics[width=0.6\textwidth]{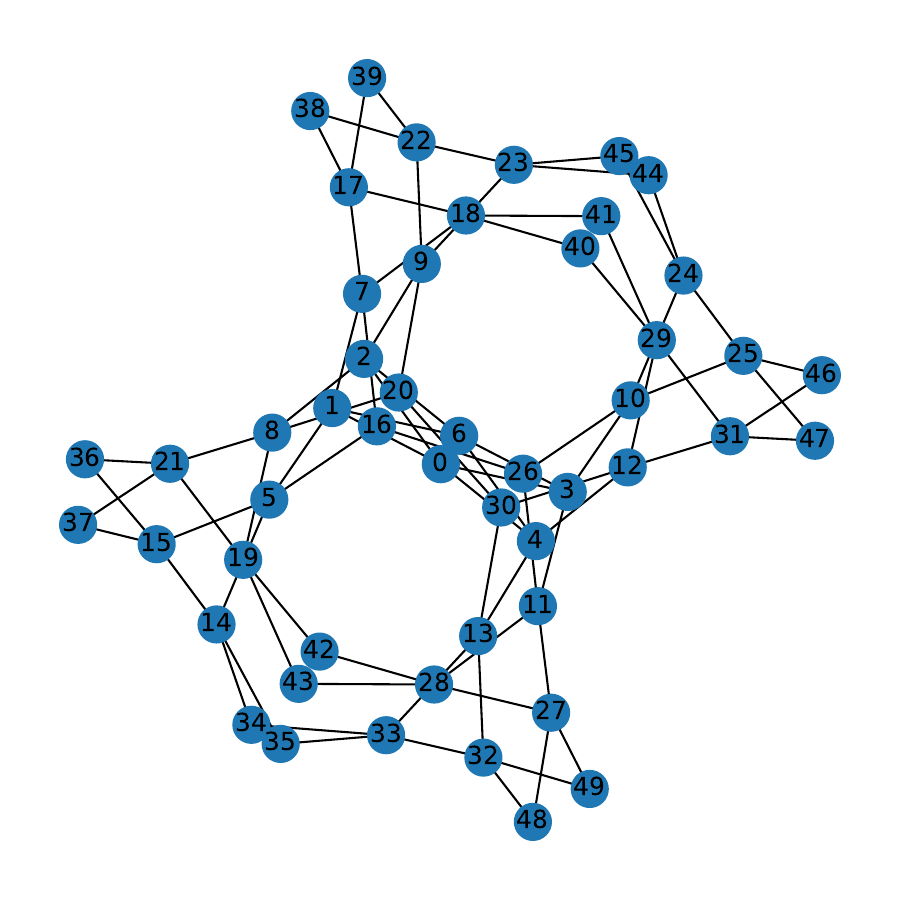}
    	\caption{A22 (finite-mutation)}\label{fig:A22EG}
    \end{subfigure} \\
	\begin{subfigure}{0.45\textwidth}
    	\centering
    	\includegraphics[width=0.6\textwidth]{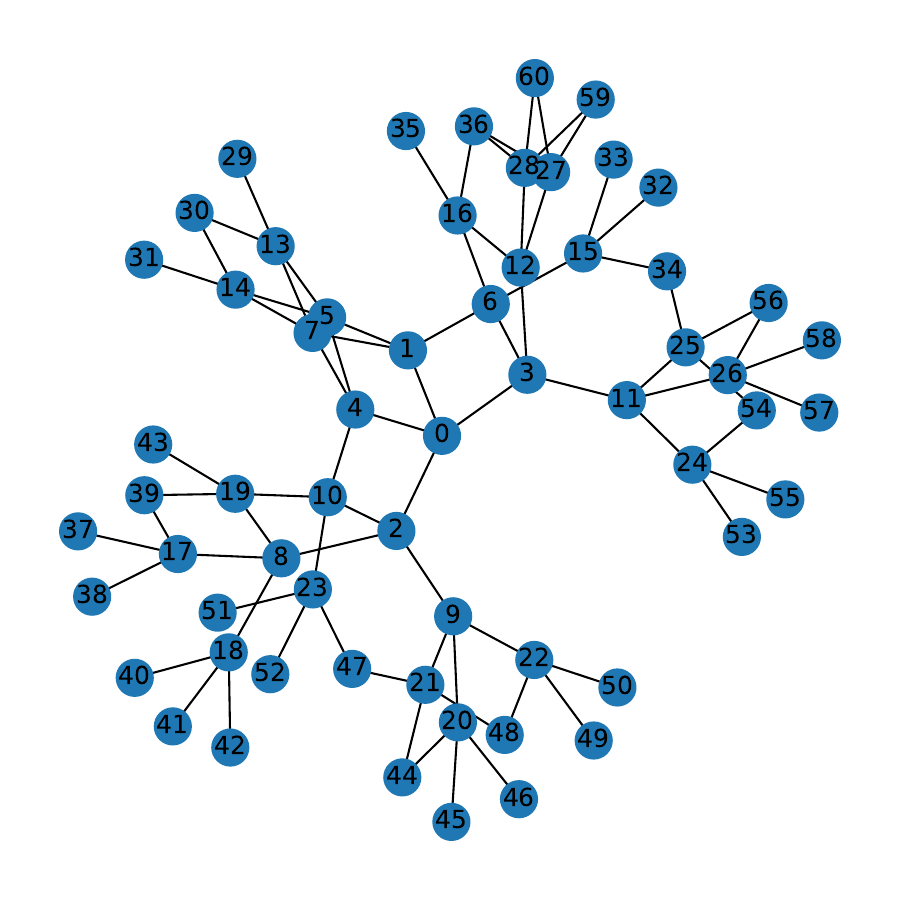}
    	\caption{$\mathcal{I}1$ (infinite)}\label{fig:inf1EG}
	\end{subfigure} 
    \begin{subfigure}{0.45\textwidth}
    	\centering
    	\includegraphics[width=0.6\textwidth]{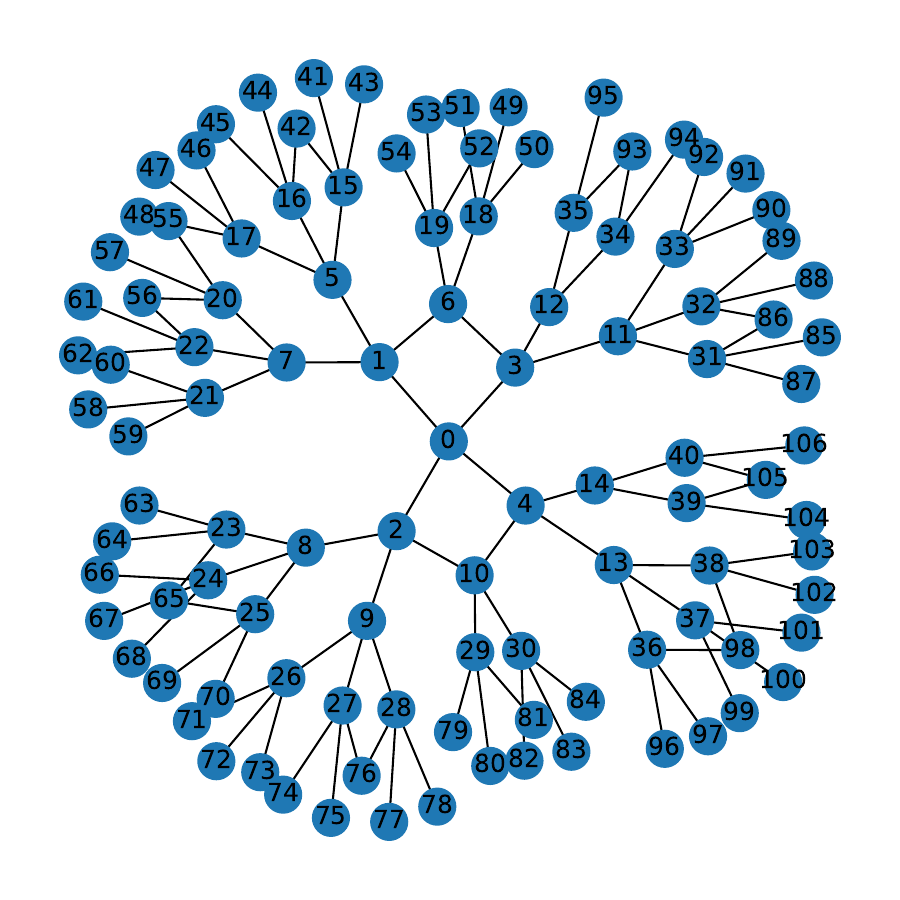}
    	\caption{$\mathcal{I}2$ (infinite)}\label{fig:inf2EG}
    \end{subfigure} \\
\caption{The quiver exchange graphs generated to depth 4 for each of the considered cluster algebras. Types are labelled, noting that finite type and finite-mutation type both have finite quiver exchange graphs. Vertices are labelled in the order they are generated starting from the initial quiver `0'.}\label{fig:QEGs}
\end{figure}

For comparison with the seed exchange graphs at depth for shown in Figure \ref{fig:EGs}, the quiver exchange graphs for the considered cluster algebras are shown plotted to depth 4 in Figure \ref{fig:QEGs}.

\subsubsection{All Types to Depth 4}
\begin{figure}[!tb]
    \centering
    \includegraphics[width=0.6\textwidth]{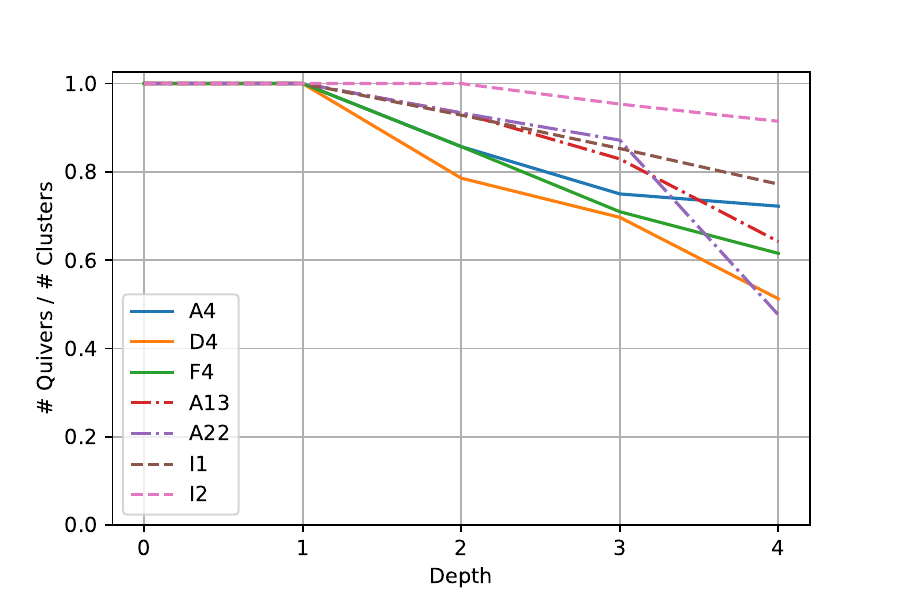}
    \caption{The relative number of quivers to clusters in the respective exchange graphs as depth varies up to depth 4 for each of the considered cluster algebras, labelled by their respective initial seeds. Each type is depicted with a different linestyle.}
    \label{fig:NumberQEGtoEG}
\end{figure}

Any embedding fundamentally depends on the number of vertices in each graph.
Therefore the number of vertices in each algebra's quiver exchange graph relative to the number in the seed exchange graph is plotted in Figure \ref{fig:NumberQEGtoEG}.
The figure shows that up to depth one there is an isomorphism between each considered algebra's quiver/seed exchange graphs, since all initial mutations will change the cluster; they hence here also always change the exchange matrix such that there are 5 distinct seeds (with distinct quivers) up to depth 1 in each case.
Since all vertices are connected, the mutated vertex will always have an orientation flip of its connected arrows, hence always changing the exchange matrix by this at least.

Beyond depth 1 the ratios drop off below 1, with the pure infinite cases dropping off the least.
These infinite cases therefore have far more distinct quivers, as may be expected since different quivers lead to different mutation processes and are hence more likely to produce new cluster variables, a defining feature of infinite type.

We can see that for all algebras there are fewer vertices and a higher density (despite $\mathcal{I}_2$ being very similar); this represents the expected behaviour from multiple seeds with different clusters having the same quivers.
Importantly, whereas where clusters are present no triangle 3-cycles can occur, when considering only quivers these triangles are possible, shown by a non-zero triangle clustering coefficient for A22.
Additionally there are more 4-cycles shown by consistently higher square clustering coefficients, which are likely used to provide alternative routes and cause the consistently smaller Wiener indices.
It is worth emphasising also that D4 has a considerably higher square clustering coefficient.

Most interestingly about the centrality analysis is that for the $\mathcal{I}_1$ algebra, the centre is no longer the initial seed!
This indicates that the extra quiver identification when the cluster information is omitted is not symmetric about the centre, and this behaviour is reflected in some of the other algebras too where the smallest difference to depth 1 is usually lower without the cluster information.

For comparison the same network analysis methods were applied to the considered algebra's quiver exchange graphs up to depth 4, as shown in Table \ref{tab:QEGanalysis}.

\begin{table}[!tb]
\centering
\addtolength{\leftskip} {-2cm}
\addtolength{\rightskip}{-2cm}
\footnotesize
\begin{tabular}{|c|cccccc|}
\hline
\multirow{2}{*}{\begin{tabular}[c]{@{}c@{}}\\ Cluster\\ Algebra\end{tabular}} & \multicolumn{6}{c|}{Quiver Exchange Graph Analysis (depth 4)} \\ \cline{2-7} 
& \multicolumn{1}{c|}{\begin{tabular}[c]{@{}c@{}}Number\\ of Vertices\end{tabular}} & \multicolumn{1}{c|}{Density} & \multicolumn{1}{c|}{\begin{tabular}[c]{@{}c@{}}Clustering\\ (tri, squ)\end{tabular}} & \multicolumn{1}{c|}{\begin{tabular}[c]{@{}c@{}}Wiener Index\\ (full, norm)\end{tabular}} & \multicolumn{1}{c|}{\begin{tabular}[c]{@{}c@{}}Centrality\\ (centre, diff)\end{tabular}} & \begin{tabular}[c]{@{}c@{}}Min cycle basis\\ ([length, freq])\end{tabular} \\ \hline
A4 & \multicolumn{1}{c|}{52} & \multicolumn{1}{c|}{0.048} & \multicolumn{1}{c|}{(0, 0.066)} & \multicolumn{1}{c|}{(6870, 5.18)} & \multicolumn{1}{c|}{(0, 0.036)} & [4,13] \\ \hline
D4 & \multicolumn{1}{c|}{41} & \multicolumn{1}{c|}{0.071} & \multicolumn{1}{c|}{(0, 0.251)} & \multicolumn{1}{c|}{(3463, 4.22)} & \multicolumn{1}{c|}{(0, 0.001)} & [4,15], [7,3] \\ \hline
F4 & \multicolumn{1}{c|}{40}  & \multicolumn{1}{c|}{0.072} & \multicolumn{1}{c|}{(0, 0.098)} & \multicolumn{1}{c|}{(3334, 4.27)} & \multicolumn{1}{c|}{(0, 0.030)} & [4,14], [6,2], [8,1] \\ \hline
A13 & \multicolumn{1}{c|}{70} & \multicolumn{1}{c|}{0.036} & \multicolumn{1}{c|}{(0, 0.041)} & \multicolumn{1}{c|}{(12826, 5.31)} & \multicolumn{1}{c|}{(0, 0.020)} & [4,9], [6,8] \\ \hline
A22 & \multicolumn{1}{c|}{50} & \multicolumn{1}{c|}{0.067} & \multicolumn{1}{c|}{(0.080, 0.108)} & \multicolumn{1}{c|}{(4780, 3.90)} & \multicolumn{1}{c|}{(0, 0.029)} & [3,8], [4,15], [7,2], [8,8] \\ \hline
$\mathcal{I}1$ & \multicolumn{1}{c|}{61} & \multicolumn{1}{c|}{0.044} & \multicolumn{1}{c|}{(0, 0.134)} & \multicolumn{1}{c|}{(9456, 5.17)} & \multicolumn{1}{c|}{(1, -)} & [4,18], [6,2] \\ \hline
$\mathcal{I}2$ & \multicolumn{1}{c|}{107} & \multicolumn{1}{c|}{0.020} & \multicolumn{1}{c|}{(0, 0.040)} & \multicolumn{1}{c|}{(33900, 5.98)} & \multicolumn{1}{c|}{(0, 0.061)} & [4,10] \\ \hline
\end{tabular}
\caption{Network analysis of the quiver exchange graphs generated to depth 4 for the 7 cluster algebras considered, labelled by their respective initial seeds. The analysis lists: the number of vertices in the EG; the density of the EG; the triangle and square average clustering coefficients; the Wiener index (both full form and normalised form); the eigenvector centrality analysis listing the central vertex and then the size of the smallest difference in centrality from the initial seed ``0'' to the clusters at depth 1 (when the initial seed is the centre); and finally the information on the minimum cycle basis showing the length of the basis cycles and the frequency of those lengths in the basis.}\label{tab:QEGanalysis}
\end{table}

Finally, the minimum cycle bases have quite non-trivial changes. 
Most notably, the cycles are no longer always even, as there is extra redundancy where mutating quivers can produce another quiver in the same depth.
In some cases the cycle bases have more cycles (D4, A13, A22, $\mathcal{I}_1$), and sometimes fewer (A4, F4, $\mathcal{I}_2$).
This behaviour is curious and highlights the subtleties of the quiver exchange graph's embedding in the seed exchange graph.
To probe this further we wish to lose the limiting behaviour of the depth truncation, and for that we return to the finite type generalised associahedra.

\subsubsection{Quiver Generalised Associahedra}
After first generating the quiver generalised associahedra, our first clue about the embedding structure comes from the number of vertices.
Whereas when identifying by the permutation-equivalence there isn't a consistent pattern, 
 when instead one does  not identify in this way, some beautiful structure emerges.

To illustrate this we show the number of vertices in the seed and quiver generalised associahedra, as well as the ratios between these numbers, for the rank 4 algebras considered here in Table \ref{tab:QEGEGratios}.

\begin{table}[!h]
\centering
\begin{tabular}{|c|c|c|c|c|c|}
\hline
Cluster Algebra          & A4   & B4  & C4  & D4   & F4  \\ \hline
Number of Quivers       & 144   & 84  & 84  & 50   & 60  \\ \hline
Number of Seeds & 1008 & 420 & 420 & 1200 & 420 \\ \hline
Ratio            & 7    & 5   & 5   & 24   & 7   \\ \hline
\end{tabular}
\caption{The number of vertices in the quiver exchange graph and seed exchange graph for the rank 4 finite type algebras considered. Their ratios are also listed, all taking integer values.}\label{tab:QEGEGratios}
\end{table}





%

A priori, one may not expect all these ratios to be integer.
When each seed in the seed generalised associahedra has its cluster information removed to leave just the quiver, there are additional identifications to be made amongst vertices where seeds with different clusters have the same quiver. However, there is no requirement for all the quivers to occur the same number of times (which is what we see and leads to the concurrent identification of each set of the `ratio' number of quivers and hence an integer ratio overall). 

These ratios take somewhat surprising and perhaps unintuitive numbers, with no obvious foundation in the Dynkin construction.
To further probe this behaviour, these ratios were also computed for the finite type algebras (arising from the 4 Dynkin series) for all ranks up to rank 5 in Table \ref{tab:ratioswithrank}, and show signs of further extraordinary structure.
Conjectured natural continuations of these observed ratios are provided for higher ranks also, noting that $D_4$ appears to be anomalous in its series, likely related to triality of its initial quiver.
\begin{table}[!t]
\centering
\begin{tabular}{|c|cccccc|}
\hline
\multirow{2}{*}{\begin{tabular}[c]{@{}c@{}}Cluster\\ Algebra\end{tabular}} & \multicolumn{6}{c|}{Rank} \\ \cline{2-7} 
& \multicolumn{1}{c|}{1} & \multicolumn{1}{c|}{2} & \multicolumn{1}{c|}{3} & \multicolumn{1}{c|}{4}  & \multicolumn{1}{c|}{5}   & $r \geq 6$      \\ \hline
$A_r$                                                                         & \multicolumn{1}{c|}{2} & \multicolumn{1}{c|}{5} & \multicolumn{1}{c|}{6} & \multicolumn{1}{c|}{7}  & \multicolumn{1}{c|}{8} & $r+3$    \\ \hline
$B_r$                                                                         & \multicolumn{1}{c|}{-} & \multicolumn{1}{c|}{3} & \multicolumn{1}{c|}{4} & \multicolumn{1}{c|}{5}  & \multicolumn{1}{c|}{6}   & $r+1$    \\ \hline
$C_r$                                                                         & \multicolumn{1}{c|}{-} & \multicolumn{1}{c|}{3} & \multicolumn{1}{c|}{4} & \multicolumn{1}{c|}{5}  & \multicolumn{1}{c|}{6}   & $r+1$    \\ \hline
$D_r$                                                                         & \multicolumn{1}{c|}{-} & \multicolumn{1}{c|}{4} & \multicolumn{1}{c|}{6} & \multicolumn{1}{c|}{24} & \multicolumn{1}{c|}{10} & $2r$ \\ \hline
\end{tabular}
\caption{The ratios between number of seeds in the seed exchange graph and number of quivers in the quiver exchange graph for the finite type cluster algebras, not applying permutation equivalence between seeds/quivers. Conjectured relationships are shown for higher ranks, $r$, beyond feasible computation at present.}\label{tab:ratioswithrank}
\end{table}
For completeness, we provide the $G_2$ ratio: 4, reiterate that the $F_4$ ratio is 7, and note that the remaining $E_6, E_7, E_8$ exceptional cases are too high a rank to be feasibly computed with current resources (leaving this to future work).

We finally reemphasise that if one identifies via the permutation equivalence some of this alluring behaviour is lost. 
If one does identify in this way, some integer ratios do occur sporadically for the $A_r$ and $D_r$ series ($A_4:7, D_5:7$) beyond the trivial rank 2 cases where there is only one quiver.
We believe this to be a probabilistic artefact occasionally carried over through the permutation identification, where there are few factors to choose from when dividing the number of seeds in these smaller rank algebras under permutation equivalence identification, so maintaining the integer ratio is more likely.
The $F_4$ and $G_2$ ratios remain 7 and 4 respectively after the permutation equivalence identification; and for the $B_r/C_r$ series all the ratios stay the same!

A final notable comment is that these conjectured ratios in Table \ref{tab:ratioswithrank} bare resemblance to selected geometric realisations of cluster algebras discussed in \cite{CA_2}.
In this, a $A_r$ cluster algebra may be realised as triangulations of a $(r+3)$-gon, $B_r$ and $C_r$ cluster algebras as arcs in an orbifolded $2(r+1)$-gon, and a $D_r$ cluster algebra as coloured diagonals in a $(2r)$-gon.
Perhaps these special realisations work due to some natural relationship to the quiver-seed exchange graph embeddings, which may in turn help shed light on the origin of these ratio values.

\paragraph{Cycle Embedding}
Already the vertex embedding of the quiver exchange graph into the seed exchange graph reveals an intriguing structure.
As motivated when examining all algebra types at depths up to 4, the cycle space also acts as a foundation of the graph structure and an important tool for analysing these embeddings.

Let us consider an $s$-cycle of quivers in a quiver exchange graph, i.e. there is a sequence of $s$ mutations connecting each quiver vertex to the next, eventually reproducing the quiver one started with.
Then considering the same algebra's seed exchange graph and taking the subgraph of seeds which have the same quivers as in the $s$-cycle produces a subgraph built out of $q$ $t$-cycles, where $t=ps$ and $p,q,s,t \in \mathbb{Z}_+$.
We call $p$ the scale factor, as it describes how the size of the cycle scales, and $q$ the copy factor, as it dictates how many copies of the cycle are produced.

Looking at all the algebras we consider, the values of $p$ and $q$ change depending on the cycle in the quiver exchange graph considered.
To best illustrate this we again focus on the generalised associahedra of the finite type rank 4 algebras, where the entire subgraph corresponding to the chosen cycle can be computed, and the value of $pq$ remains constant, equalling the respective `ratio' values calculated previously.

The embedding of each quiver exchange graph cycle can be considered by taking the subgraph of all seeds in the seed exchange graph that have a quiver from the cycle under consideration.
For all the rank 3 and 4 finite type cluster algebras the distribution of $p$ and $q$ values for the cycles in each algebra's minimum cycle basis are given in Table \ref{tab:qmcb_embedding}.
The minimum cycle basis was used here as it is a sensible set of independent cycles of different sizes to probe the embedding behaviour.

\begin{table}[!tb]
\centering
\addtolength{\leftskip} {-2cm}
\addtolength{\rightskip}{-2cm}
\footnotesize
\begin{tabular}{|c|c|c|c|}
\hline
Cluster Algebra & \begin{tabular}[c]{@{}c@{}}QEG MCB\\ $[[$len,freq$]]$\end{tabular} & \begin{tabular}[c]{@{}c@{}}Cycle scale factor $p$\\ $[[$value,freq$]]$\end{tabular} & \begin{tabular}[c]{@{}c@{}}Cycle copy factor $q$\\ $[[$value,freq$]]$\end{tabular} \\ \hline
A4 & [[4,108],[6,8],[10,29]] & [[1,90],[7,55]] & [[1,55],[7,90]] \\ \hline
B4=C4 & [[4,60],[6,15],[8,6],[10,4]] & [[1,49],[5,36]] & [[1,36],[5,49]] \\ \hline
D4  & [[4,33],[7,12],[8,6]] & [[1,15],[2,2],[4,25],[6,6],[12,3]] & [[2,3],[4,6],[6,25],[12,2],[24,15]] \\ \hline
F4 & [[4,42],[6,14],[8,4],[10,1]] & [[1,31],[7,30]] & [[1,30],[7,31]] \\ \hline
A3=D3 & [[3,6],[8,2]] & [[3,2],[6,6]] & [[1,6],[2,2]] \\ \hline
B3=C3 & [[3,4],[5,2]] & [[4,6]] & [[1,6]] \\ \hline
\end{tabular}
\caption{The embedding of the quiver exchange graph minimum cycle basis (QEG MCB) into the respective seed exchange graph for the rank 3 \& 4 finite type cluster algebras. The embedding information is listed as the $p$ \& $q$ values and frequencies that dictate how each cycle scales in size and copies respectively.}\label{tab:qmcb_embedding}
\end{table}

As can be seen, the $p$ and $q$ values are not constant for each algebra, but matching up the frequencies (and computationally confirmed explicitly) shows that the $pq$ ratio is always constant at the value listed in Table \ref{tab:ratioswithrank}.
Due to the requirement that all cycles in the seed exchange graph are even, if the quiver exchange graph cycle being embedded is odd then $p$ has to be even, which is well exemplified with the rank 3 cases where for A3 the 6 3-cycles have $p=6$, and all the B3 cycles have $p=4$ (i.e. $p=1$ cannot occur).
Note this breaks the symmetry of the ratio being split into its factors where $p$ and $q$ can take either factor's value, as seen for the even cycles in the table.

\begin{figure}[!tb]
	\centering
	\begin{subfigure}{0.45\textwidth}
    	\centering
    	\includegraphics[width=0.6\textwidth]{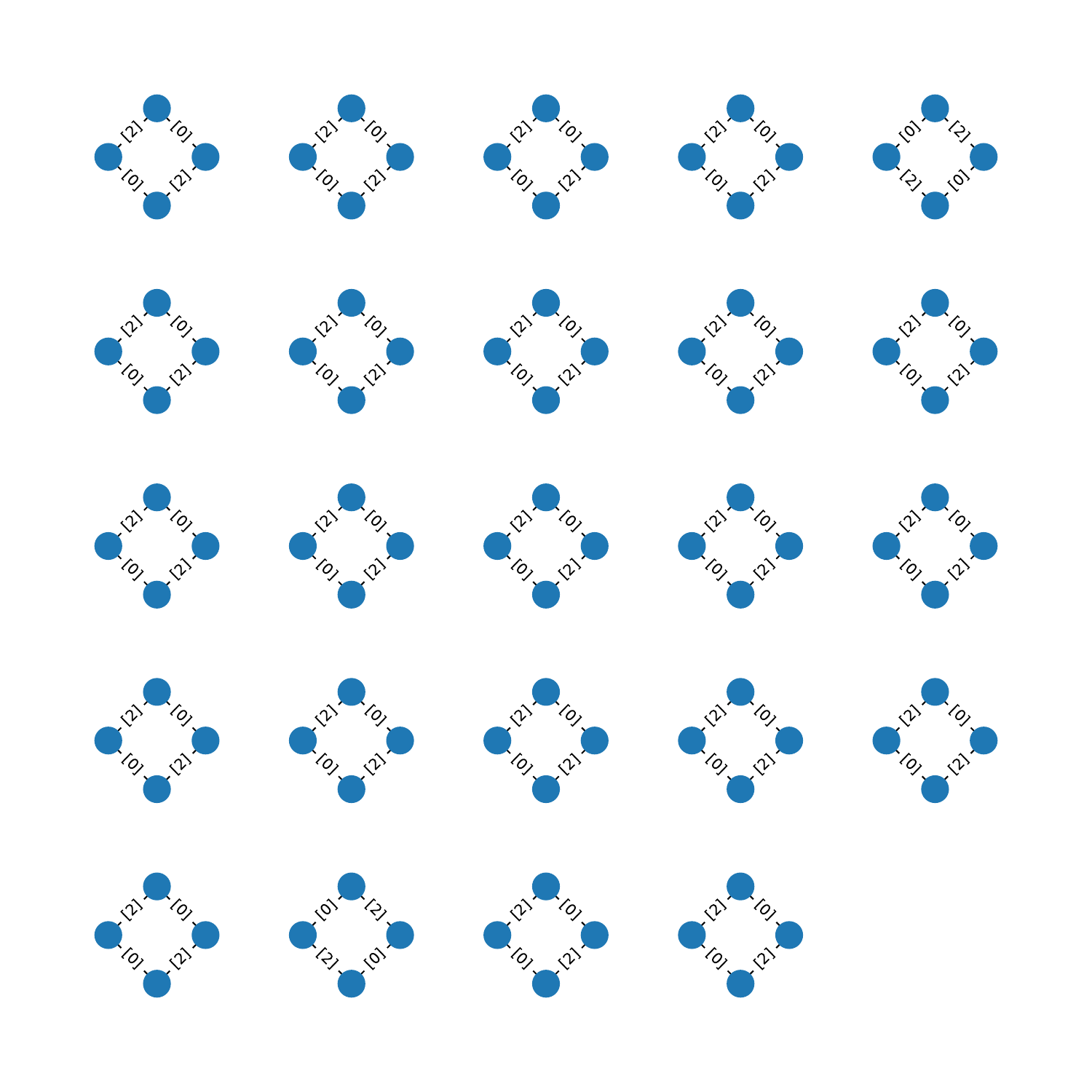}
    	\caption{QEG 4-cycle ($p=1,q=24$)}\label{fig:D4_4cyclecommute}
	\end{subfigure}
    \begin{subfigure}{0.45\textwidth}
    	\centering
    	\includegraphics[width=0.6\textwidth]{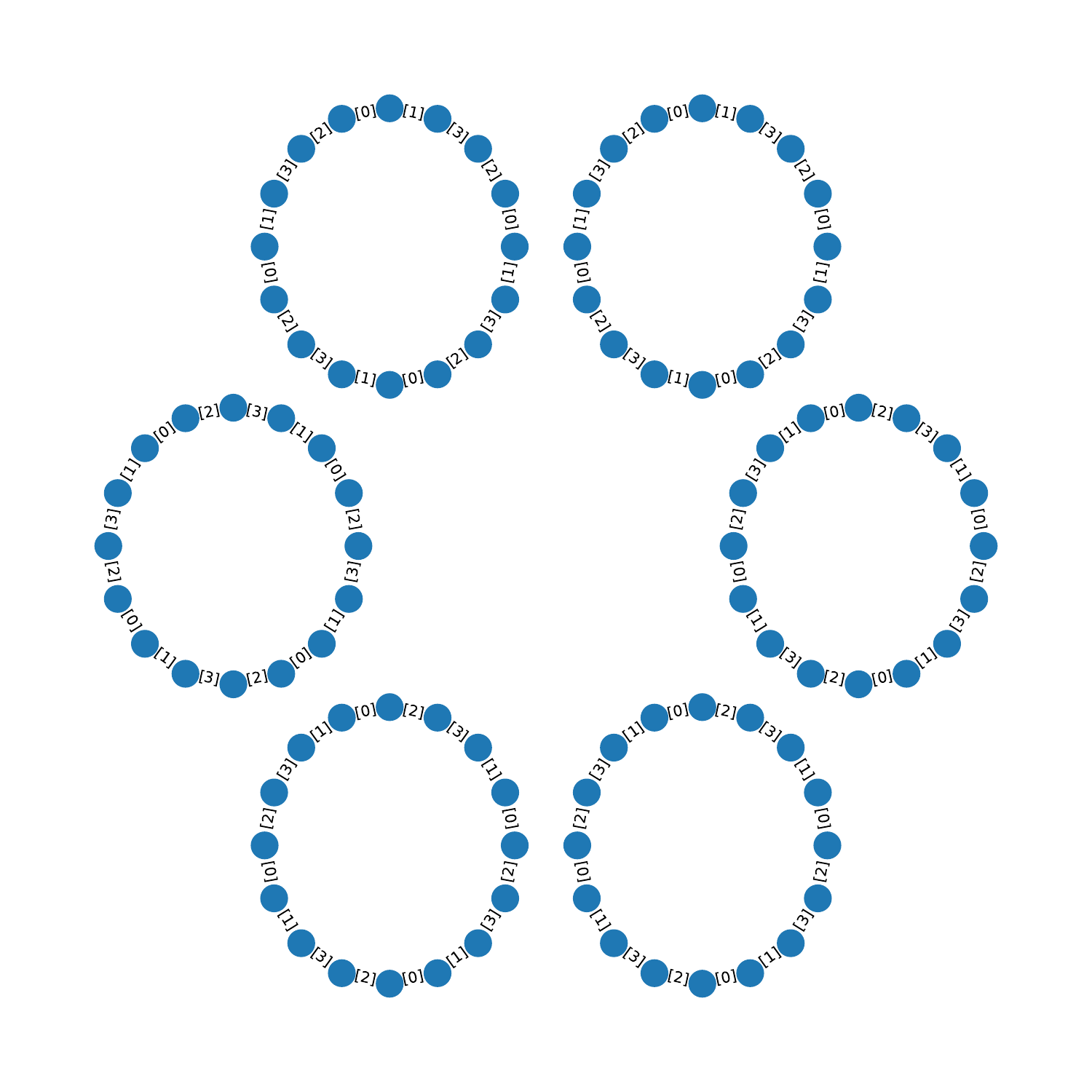}
    	\caption{QEG 4-cycle ($p=4,q=6$)}\label{fig:D4_4cyclenoncommute}
    \end{subfigure}  \\
    \begin{subfigure}{0.45\textwidth}
    	\centering
    	\includegraphics[width=0.6\textwidth]{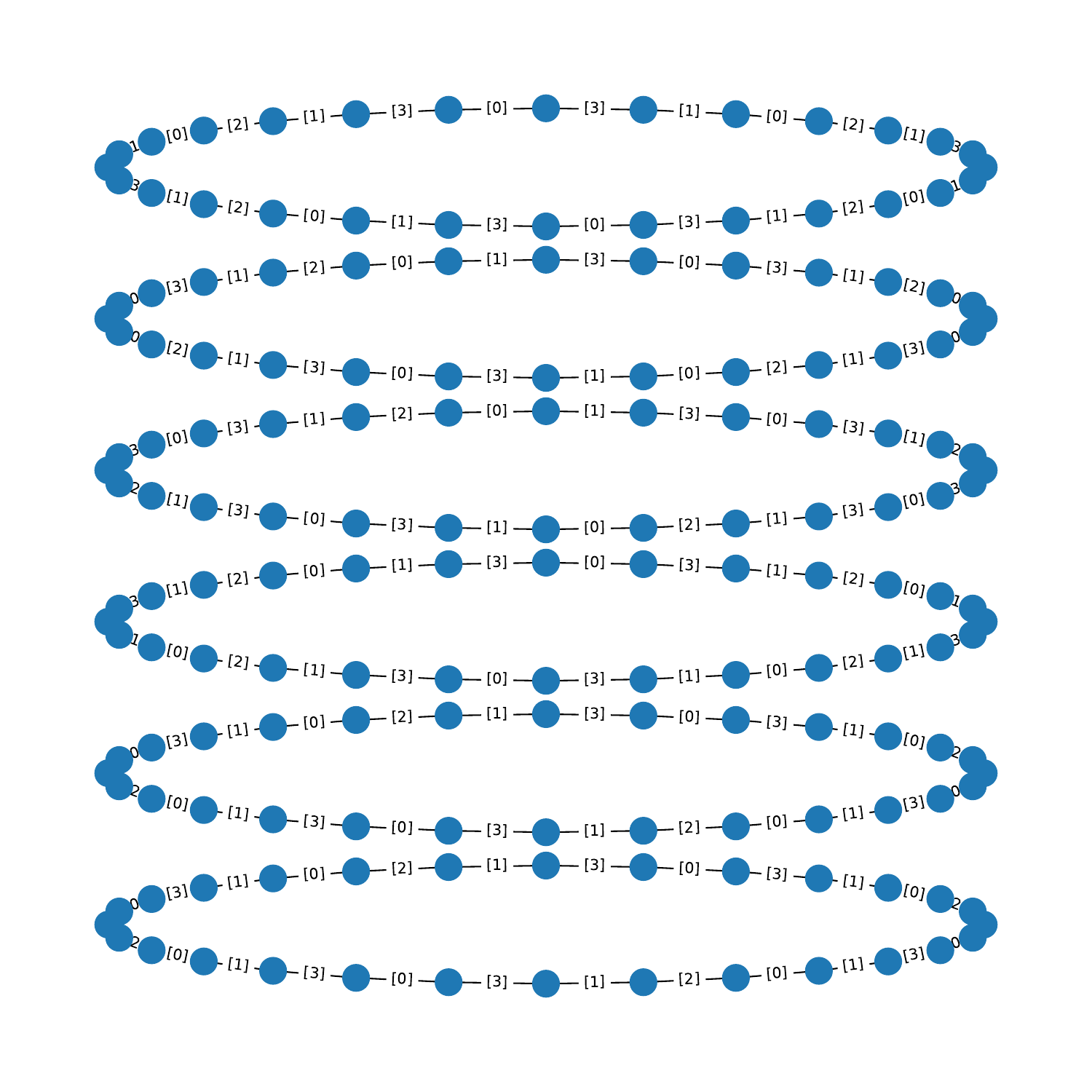}
    	\caption{QEG 7-cycle ($p=4,q=6$)}\label{fig:D4_7cycle}
    \end{subfigure}
\caption{The seed exchange graph embedding of 3 selected cycles from the quiver exchange graph (QEG) for the D4 cluster algebra. (a) \& (b) are 4-cycles, which when embedded become 24 4-cycles and 6 16-cycles respectively; whilst (c) is a 7-cycle which becomes 6 28-cycles. Each embedded cycle subgraph has edge features indicating the respective mutating quiver vertex, which connects the seeds, forming multiples of the quiver 4 and 7-cycles respectively.}\label{fig:D4_cycleembed}
\end{figure}

To provide some explicit example of the embedding we focus on the D4 algebra.
The embedding of 3 different quiver exchange graph cycles in the seed exchange graph are given in Figure \ref{fig:D4_cycleembed}.
As can be seen, the two 4-cycles have different $(p,q)$ values, and the larger 7-cycle has an even $p$ value.

Of particular note is the 4-cycle with seed exchange graph embedding shown in Figure \ref{fig:D4_4cyclecommute}. 
This quiver exchange graph 4-cycle comes from commuting action of mutation on two different vertices, which in the quiver happens when the vertices mutated about are not connected.
When any 4-cycle (in any cluster algebra) is built from commuting action on disconnected vertices, the respective embedding in the seed exchange graph always has $p=1$, such that the mutation remains a commuting relation with the cluster information added for all clusters with those quivers.
Since the two mutated quiver vertices are not connected, then no edge connecting them is introduced, and in the cluster each new variable from a mutation will not include the variable from the unconnected vertex.
This does not hold for all generic 4-cycles, as exemplified in Figure \ref{fig:D4_4cyclenoncommute} where the quiver 4-cycle comes from mutation on different vertices and has $p>1$; therefore only when the cycle is from commuting action will $p=1$.

In addition, whilst the finite type cases have constant $pq$, putting an upper-bound on the $p$ and $q$ values for all cycles in the quiver exchange graph, in the finite-mutation but not finite and infinite types there is no such upper bound.
In fact, particularly for the finite-mutation but not finite cases, since there are finitely many quivers but infinitely many seeds, for some cycles either (or both) of the $p$ or $q$ values must be infinite.

From examining these finite type cluster algebras we have seen that the $pq$ value remains constant for all cycles; where the quiver exchange graph cycle is odd then $p$ must be even, and where a quiver exchange graph 4-cycle is from commuting action of mutation then $p=1$.
How the $p$ and $q$ values are determined more generally for each cycle in each algebra is left open for future exploratory work.

\section{Machine Learning}\label{ML}
Given two seeds it is unclear as to whether a sequence of mutations exists which connects them, i.e. whether they belong to the same cluster algebra.
Beyond simple checks for the mutation type (a necessary but not sufficient condition), brute force computation of all mutations is the usual method for checking this equivalence.
However, brute force mutation is extremely computationally expensive, as particularly emphasised by the infeasibility of computing the exchange graph for $\mathcal{I}2$ beyond depth 4.

To speed up this process of checking mutation equivalence, here machine learning methods are applied to this problem, with the idea they may be able to find invariants in the variables' Laurent polynomial structure under mutation and use it for this speeding up of equivalence computation.

In this work, dense feed-forward neural networks were used -- as the prototypical supervised machine learning classifier architecture \cite{scikit-learn}.
The neural networks had 3 layers of 256 neurons with ReLU activation, learning the binary cross-entropy loss with the Adam optimiser \cite{kingma2017adam} on these binary classification problems.
Learning performance was measured with the accuracy and Matthew's correlation coefficient (MCC) metrics, to determine the proportion of correctly classified seeds (where MCC also accounts for data bias).
Furthermore, 5-fold cross-validation was performed, training 5 independent neural networks on 5 different 80\% partitions of the dataset, such that the metrics could be averaged and standard error calculated to provide confidence in the results.

\paragraph{Data Representation}
In order to be used as input to a neural network the seed information needs to be represented by a tensor. Before representing a seed as a tensor, each cluster variable first needs to be represented, which requires encoding each Laurent polynomial's information.
The information is split into numerator and denominator, noting that the denominator will always be a single monomial term due to the Laurent phenomenon (even if it is trivially 1), whilst the numerator will be a \textit{positive}\footnote{All cluster variables of all algebras considered in this study have exclusively positive coefficients.} sum of monomials.
Each monomial, $x_1^\alpha x_2^\beta x_3^\gamma x_4^\delta$, is represented in a rank 4 tensor as the entry $T_{\alpha\beta\gamma\delta}$, with value equal to the monomial's coefficient.
Then the full variable is represented with two tensors (that are flattened and concatenated): one for all numerator monomials and one for the denominator monomial.

However, this representation is excessively sparse, for example with a proportion of non-zero entries $\sim 0.000008$ for the $\mathcal{I}1$ dataset.
Since this representation has an excessive redundancy in information, we turn to a sparse data representation, based on the `coo' style.
Here each monomial is represented as a 5-vector of entries $[c,\alpha,\beta,\gamma,\delta]$ for monomial with coefficient $c$.
Then all the monomial entries are concatenated with the denominator monomial at the end.
For comparison the sparsity proportion of the $\mathcal{I}1$ data with this style was $\sim 0.02$, a significant improvement for ML implementation.

After all tensors are generated for each variable in a cluster, they are flattened and concatenated across the cluster with the flattened exchange matrix also, to produce a single data vector per seed.
Note that sometimes the exchange matrix was omitted from the representation in order to examine the ML performance based on only the clusters, i.e. can the neural networks distinguish the cluster algebra from just sets of generators (or is the exchange matrix needed as well in order to learn the full generator set structure). Note, however, that the initial seed's cluster is the same for all the algebras, and there is further repetition also; this leads to some redundancy in information where the same tensor may be affiliated to multiple algebras.
This should make the learning noticeably harder without the exchange matrix information, as some data will mislead the learning.

Since different variables have different numbers of monomial terms, the length of the seeds' vectors varies substantially.
To create a consistent input length for all vectors in an investigation the tensors were post-padded with zeros such that all vectors were the same length as the longest in that investigation.
As the infinite type mutations lead to more complex Laurent polynomials with far higher degrees, the respective vectors for these higher-depth outer seeds are much longer than all others across all the algebras.
Therefore in order to stop this dilution of the information for ML comparison between non-infinite algebras, the investigations were designed to be binary classifications.

For all ML investigations seeds in the same algebra were not considered identical if they were related via a permutation of variables in the cluster. 
Due to the systematic vector generation procedure this led to different vectors for different permutations, and one may consider this process as the common practice of data augmentation on the permutation-invariant algebra seeds if preference is to consider these over those where different permutations are unique.
The data, as well as data generation scripts, are available with this work's respective  \href{https://github.com/edhirst/ClusterAlgebrasML.git}{GitHub}.

\subsection{Distinguishing Cluster Algebra Types}
The first investigation uses the above-described neural network architecture to classify seeds coming from different cluster algebras.
The algebras considered in this work amount to 3 finite type, 2 finite-mutation but not finite type, and 2 infinite type.
These numbers were selected such that binary classifications could be performed for all combinations of types, and in particular the additional F4 finite type was introduced, as its tensor data is closer in form to the infinite types (with more larger-than-unit entries, due to the initial quiver's non-simply laced edge double multiplicity).
Furthermore, as there are many infinite-type initial seeds, $\mathcal{I}1$ was also specifically chosen due to its lower quiver edge multiplicities, making it more similar to the finite type cases.

These selections of the algebras considered were all made such that the tensor representations of each algebra could not be distinguished by eye; one may look at this data directly in the respective \href{https://github.com/edhirst/ClusterAlgebrasML.git}{GitHub}.
Therefore any learning results would be non-trivial and constitute some true learning of the algebra structure.

For all pairs  selected  for binary classification, both algebras were generated to depth 4, their seeds converted to vectors, both sets of vectors shuffled together, and the 5-fold cross-validation ML performed.
Learning results for each of these investigations are provided in Table \ref{tab:ml_results}, where each investigation is repeated both with and without the exchange matrix information (removal reducing the vector length by 16 each time). 

\begin{table}[!tb]
\centering
\footnotesize
\begin{tabular}{|c|c|c|c|cccc|}
\hline
\multirow{3}{*}{Investigation} & \multirow{3}{*}{Class Sizes} & \multirow{3}{*}{\begin{tabular}[c]{@{}c@{}}Tensor \\ Length\end{tabular}} & \multirow{3}{*}{\begin{tabular}[c]{@{}c@{}}Tensor \\ Sparsity\end{tabular}} & \multicolumn{4}{c|}{ML Performance} \\ \cline{5-8} 
& & & & \multicolumn{2}{c|}{with EM} & \multicolumn{2}{c|}{no EM} \\ \cline{5-8} & & & & \multicolumn{1}{c|}{Accuracy} & \multicolumn{1}{c|}{MCC} & \multicolumn{1}{c|}{Accuracy} & MCC \\ \hline
A4 vs D4 & (72, 80) & 180 & 0.120 & \multicolumn{1}{c|}{\begin{tabular}[c]{@{}c@{}}0.867\\ $\pm$ 0.021\end{tabular}} & \multicolumn{1}{c|}{\begin{tabular}[c]{@{}c@{}}0.741\\ $\pm$ 0.036\end{tabular}}  & \multicolumn{1}{c|}{\begin{tabular}[c]{@{}c@{}}0.893\\ $\pm$ 0.026\end{tabular}} & \begin{tabular}[c]{@{}c@{}}0.788\\ $\pm$ 0.053\end{tabular} \\ \hline
A4 vs A13 & (72, 109) & 280 & 0.088 & \multicolumn{1}{c|}{\begin{tabular}[c]{@{}c@{}}0.944\\ $\pm$ 0.011\end{tabular}} & \multicolumn{1}{c|}{\begin{tabular}[c]{@{}c@{}}0.886\\ $\pm$ 0.023\end{tabular}}  & \multicolumn{1}{c|}{\begin{tabular}[c]{@{}c@{}}0.878\\ $\pm$ 0.022\end{tabular}} & \begin{tabular}[c]{@{}c@{}}0.743\\ $\pm$ 0.047\end{tabular} \\ \hline
F4 vs $\mathcal{I}$1 & (65, 79) & 2320 & 0.015 & \multicolumn{1}{c|}{\begin{tabular}[c]{@{}c@{}}0.950\\ $\pm$ 0.013\end{tabular}} & \multicolumn{1}{c|}{\begin{tabular}[c]{@{}c@{}}0.903\\ $\pm$ 0.024\end{tabular}} & \multicolumn{1}{c|}{\begin{tabular}[c]{@{}c@{}}0.936\\ $\pm$ 0.012\end{tabular}} & \begin{tabular}[c]{@{}c@{}}0.875\\ $\pm$ 0.024\end{tabular} \\ \hline
A13 vs A22 & (109, 105) & 280 & 0.091 & \multicolumn{1}{c|}{\begin{tabular}[c]{@{}c@{}}0.810\\ $\pm$ 0.028\end{tabular}} & \multicolumn{1}{c|}{\begin{tabular}[c]{@{}c@{}}0.630\\ 0.050\end{tabular}}     & \multicolumn{1}{c|}{\begin{tabular}[c]{@{}c@{}}0.810\\ $\pm$ 0.024\end{tabular}} & \begin{tabular}[c]{@{}c@{}}0.633\\ $\pm$ 0.049\end{tabular} \\ \hline
A13 vs $\mathcal{I}$1 & (109, 79) & 2320 & 0.015 & \multicolumn{1}{c|}{\begin{tabular}[c]{@{}c@{}}0.930\\ $\pm$ 0.021\end{tabular}} & \multicolumn{1}{c|}{\begin{tabular}[c]{@{}c@{}}0.855\\ $\pm$ 0.048\end{tabular}}  & \multicolumn{1}{c|}{\begin{tabular}[c]{@{}c@{}}0.914\\ $\pm$ 0.021\end{tabular}} & \begin{tabular}[c]{@{}c@{}}0.801\\ $\pm$ 0.059\end{tabular} \\ \hline
$\mathcal{I}$1 vs $\mathcal{I}$2 & (79, 117) & 94280 & 0.008  & \multicolumn{1}{c|}{\begin{tabular}[c]{@{}c@{}}0.918\\ $\pm$ 0.023\end{tabular}} & \multicolumn{1}{c|}{\begin{tabular}[c]{@{}c@{}}0.830\\ $\pm$ 0.047\end{tabular}}  & \multicolumn{1}{c|}{\begin{tabular}[c]{@{}c@{}}0.923\\ $\pm$ 0.021\end{tabular}} & \begin{tabular}[c]{@{}c@{}}0.840\\ $\pm$ 0.043\end{tabular} \\ \hline
\end{tabular}
\caption{Machine learning results for NN binary classification between clusters generated by the respective initial seeds. The sizes of each class are listed, along with the tensor length used to represent them, and the average sparsity of those tensors (proportion of non-zero entries). The investigations are carried out with and without the exchange matrix (EM) information for each cluster. The performance is measured by accuracy and MCC with 5-fold cross-validation to provide standard error confidence on the measures.}\label{tab:ml_results}
\end{table}

For each investigation the  algebras considered are labelled by their initial seeds; the respective class sizes are given, with the full tensor length into which all the seeds are embedded (based on the largest seed in that investigation).
In order to compare investigations the tensor sparsity is also given.
As can be seen from this meta-data, the data size for training is very small relative to usual ML investigations, and especially with the low proportion of non-zero terms there is little information for a neural network to learn any relationship between seeds from the same cluster algebra.

It is therefore evermore surprising that the architecture learns so well in all investigations, with accuracies and MCC scores $>90\%$ in some investigations.
It can therefore be confidently concluded that machine learning can identify structure inherent to each algebra, and learn the cluster mutation process.

Interestingly, the inclusion of the exchange matrix information only improves learning for classifications between algebras of different type (finite, finite-mutation, infinite).
This may be due to the quivers for algebras of the same type looking more similar and therefore diluting the already sparse information in the tensors, hindering the learning performance.
This is also surprising since removal of the exchange matrix information makes some seeds in different algebras identical (for example the initial seeds); hence one would expect performance to always be worse without it.
It must therefore be the case that the dilution of relevant information is a more substantial factor than occasional misleading of the learning for same type classifications.

\paragraph{Distinguishing Generalised Associahedra}
Whereas the previous investigations focused on all types with data generated to depth 4, here we learn only the finite types -- but generated to their maximum depth to include all seeds in the algebras.

Since these all have similar complexity cluster variables we embed them all in the same size tensors of length 1576 (with the exchange matrix information), and perform multiclassification between all of them (including B4 and C4 also).
The neural network architecture is the same except now cross-entropy loss must be used instead of binary cross-entropy.

The 5-fold cross-validation results give averaged performance measures:
\begin{align}
    Accuracy & = \qquad\qquad 0.989 \pm 0.003 \;,\\
    MCC & = \qquad\qquad 0.985 \pm 0.004 \;,\\
    CM & = \begin{pmatrix} 
            0.289 & 0.003 & 0.000 & 0.001 & 0.000 \\
            0.001 & 0.120 & 0.000 & 0.000 & 0.000 \\
            0.002 & 0.000 & 0.119 & 0.000 & 0.000 \\
            0.002 & 0.000 & 0.000 & 0.344 & 0.000 \\
            0.004 & 0.001 & 0.000 & 0.000 & 0.117
            \end{pmatrix}\;.
\end{align}
The learning performance is exceptionally strong, despite now demanding multiclassification from the architecture.

The averaged confusion matrix, $CM$, shows the proportion of truly class A4:B4: C4:D4:F4 (given by the row) classified into class A4:B4:C4:D4:F4 (given by the column).
Perfect learning produces a diagonal matrix, and here the learning is very close to that, with off-diagonal components two orders of magnitude smaller than the diagonal components.
The larger diagonal entries of A4 and D4 reflect their larger frequencies in the dataset.
The most frequent non-zero off-diagonals occur where the other algebras are more likely to misclassify as A4 (larger first column entries), potentially due to lower depth seeds all being exceptionally similar across the algebras which the neural network then arbitrarily assumes to all be A4.
Surprisingly, the matrix shows that the architecture can distinguish well between the B4 and C4 architectures, despite the analysis showing that they have identical generalised associahedra structure.

\subsection{Learning at Varying Depths}
In order to connect these results of learning at depth 4 and at the maximum depth for the finite-type full algebras, we examine the A4:D4 binary classification performance as depth increases from the minimum possible depth 1 (such that enough data to train \& test) up to depth 13 (where both algebras have all their seeds generated).
Here we consider the ML investigation without the exchange matrix information, as the results in the previous section suggested it was in some sense superfluous for the learning.

As the depth increases not only do the cluster variables become more complex and hence represented by longer and sparser vectors, there are more variables to train with too.
This investigation aims to probe these competing effects of longer, more complicated vectors to learn from, against the benefit of more data to learn with.
The information regarding the tensor length needed for embedding, and the respective A4 and D4 class sizes at each depth are given in Table \ref{tab:mldata_varydepth}.

\begin{figure}[!tb]
    \centering
    \includegraphics[width=0.6\textwidth]{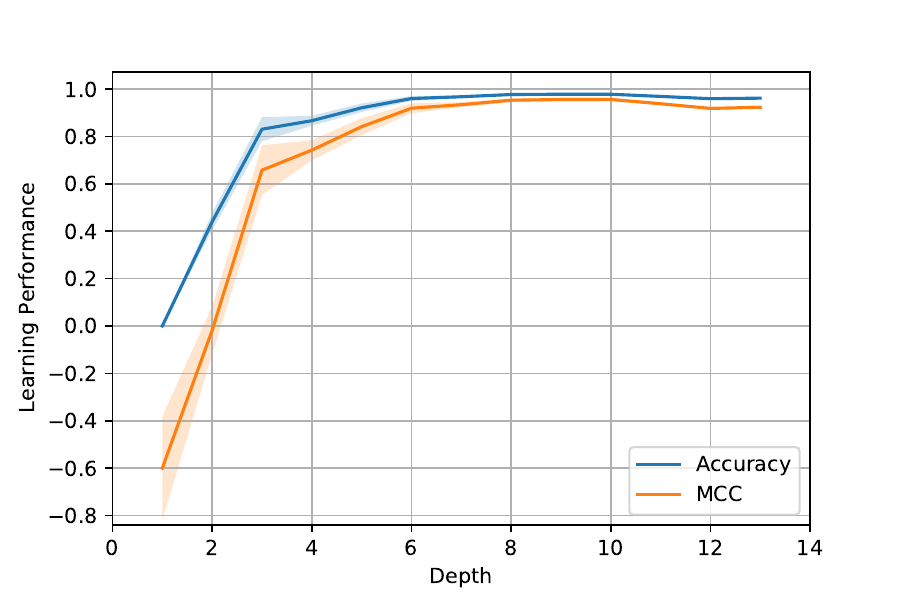}
    \caption{ML results for binary classification between A4 and D4 clusters with data generated from the respective initial seeds up to the given depths. Performance is measured with accuracy and MCC, and over the 5-fold cross-validation the measures are averaged and standard error calculated giving the shown error bounds.}
    \label{fig:A4D4_varydepth}
\end{figure}

The cross-validation learning results are this time plotted as depth varies in Figure \ref{fig:A4D4_varydepth}, with the performance measures' standard errors given as error bounds.
It can be clearly seen that as depth increases the NNs perform better in the classification, even at lower depths where the tensor length has not yet stabilised.

\begin{table}[!h]
\centering
\footnotesize
\begin{tabular}{|c|c|c|c|c|c|c|c|c|c|c|c|c|c|}
\hline
Depth          & 1 & 2 & 3 & 4 & 5 & 6 & 7 & 8 & 9 & 10 & 11 & 12 & 13  \\ \hline
Class Sizes    & \begin{tabular}[c]{@{}c@{}}5\\ 5\end{tabular} & \begin{tabular}[c]{@{}c@{}}14\\ 14\end{tabular} & \begin{tabular}[c]{@{}c@{}}32\\ 33\end{tabular} & \begin{tabular}[c]{@{}c@{}}72\\ 80\end{tabular} & \begin{tabular}[c]{@{}c@{}}151\\ 180\end{tabular} & \begin{tabular}[c]{@{}c@{}}283\\ 372\end{tabular} & \begin{tabular}[c]{@{}c@{}}462\\ 658\end{tabular} & \begin{tabular}[c]{@{}c@{}}653\\ 928\end{tabular} & \begin{tabular}[c]{@{}c@{}}815\\ 1091\end{tabular} & \begin{tabular}[c]{@{}c@{}}927\\ 1167\end{tabular} & \begin{tabular}[c]{@{}c@{}}988\\ 1195\end{tabular} & \begin{tabular}[c]{@{}c@{}}1007\\ 1200\end{tabular} & \begin{tabular}[c]{@{}c@{}}1008\\ 1200\end{tabular}\\ \hline
Lengths & 76 & 96 & 136 & 196 & 196 & 196 & 196 & 196 & 196 & 196 & 196 & 196 & 196  \\ \hline
\end{tabular}
\caption{Data information for the binary classification between A4 clusters and D4 clusters generated for depths 1-13 (such that all clusters were generated). The class sizes for A4 are shown above those for D4 respectively, as well as the lengths of the flattened tensors that the clusters are embedded in.}\label{tab:mldata_varydepth}
\end{table}

\subsection{Identifying Cluster Algebra Seeds in Random Tensors}
Whereas preceding investigations used ML to differentiate which cluster algebra a specific seed generates, it is also interesting to see if neural networks can identify tensor representations which represent sensible seeds altogether.
An easy example would be representations whose exchange matrices are not skew-symmetrizable, or have diagonal elements; but beyond these there are many ways the encoding could define nonsensical seeds (especially if adding non-zeroes deep into the padding).

In order to explore this learning we first generate suitable fake data. 
To ensure the fake data is representative, and not trivially distinguishable by eye, the true seed data for each cluster variable is analysed.
Each algebra is generated on its own to depth 4, reformulated as a vector, and padded to the maximum vector length for that algebra. 
The set of vectors for each algebra is then assessed to give a discrete distribution of frequencies of all integer entries that occur in all seed vectors across the algebra.

Then as many fake vectors are generated as there are true vectors, generated to be the same lengths as the true, with each entry's value drawn from the respective discrete distribution of possible values. 
All the fake vectors were checked to not overlap with the true vectors, despite a highly improbable chance of this occurring.
For each  algebra considered the datasets of true and fake vectors were shuffled and binary classification performed, with results shown in Table \ref{tab:ml_fake}.

\begin{table}[!tb]
\centering
\footnotesize
\begin{tabular}{|c|c|c|c|c|c|c|c|} 
\hline
\multirow{3}{*}{\begin{tabular}[c]{@{}c@{}}\\Performance\\ Measure\end{tabular}} & \multicolumn{7}{c|}{Cluster Algebra}  \\ 
\cline{2-8} & A4 & D4 & F4 & A13 & A22 & $\mathcal{I}$1 & $\mathcal{I}$2 \\
& \multicolumn{1}{c|}{(196)} & \multicolumn{1}{c|}{(136)} & \multicolumn{1}{c|}{(336)} & \multicolumn{1}{c|}{(296)} & \multicolumn{1}{c|}{(176)} & \multicolumn{1}{c|}{(2336)} & \multicolumn{1}{c|}{(94296)} \\ 
\hline
Accuracy & \begin{tabular}[c]{@{}c@{}}1.000\\ $\pm$ 0.000\end{tabular} & \begin{tabular}[c]{@{}c@{}}1.000\\ $\pm$ 0.000\end{tabular} & \begin{tabular}[c]{@{}c@{}}1.000\\ $\pm$ 0.000\end{tabular} & \begin{tabular}[c]{@{}c@{}}1.000\\ $\pm$ 0.000\end{tabular} & \begin{tabular}[c]{@{}c@{}}1.000\\ $\pm$ 0.000\end{tabular} & \begin{tabular}[c]{@{}c@{}}0.819\\ $\pm$ 0.0488\end{tabular} & \begin{tabular}[c]{@{}c@{}}0.800\\ $\pm$ 0.028\end{tabular}  \\ 
\hline
MCC & \begin{tabular}[c]{@{}c@{}}1.000\\ $\pm$ 0.000\end{tabular} & \begin{tabular}[c]{@{}c@{}}1.000\\ $\pm$ 0.000\end{tabular} & \begin{tabular}[c]{@{}c@{}}1.000\\ $\pm$ 0.000\end{tabular} & \begin{tabular}[c]{@{}c@{}}1.000\\ $\pm$ 0.000\end{tabular} & \begin{tabular}[c]{@{}c@{}}1.000\\ $\pm$ 0.000\end{tabular} & \begin{tabular}[c]{@{}c@{}}0.671\\ $\pm$ 0.078\end{tabular}  & \begin{tabular}[c]{@{}c@{}}0.640\\ $\pm$ 0.044\end{tabular}  \\
\hline
\end{tabular}
\caption{Binary classification results for differentiating tensors representing cluster algebras generated to depth 4 from the respective listed initial seeds, against fake tensors generated to mimic them. The respective tensor lengths are listed beneath the initial seeds in brackets. Performance is measured with accuracy and MCC across the 5-fold cross-validation runs.}\label{tab:ml_fake}
\end{table}

The results show perfect classification for all except the infinite types.
This indicates that the neural networks can learn some non-trivial structure in the finite and finite-mutation types which it can use to effectively differentiate from fake data.
However, for the infinite types the poorer performance suggests that the tensor structure is perhaps more erratic and hence harder to differentiate from the  simple random uniform model for its fake data.
This infinite data may be expected to span a larger proportion of the possible tensors generated since there are infinitely many of them in the algebra, and hence it may also be the case that these fake tensors are related to seeds at higher depths.

\section{Summary \& Outlook}\label{summary}
Network analysis methods uncovered patterns in the cluster algebra exchange graphs unique to each type.
In particular, a symmetric behaviour for quiver exchange graph embedding in the seed exchange graph showed constant integer ratios between respective numbers of graph vertices, which we conjecture for Dynkin types of any rank.
This behaviour is made manifest by omitting the permutation-equivalence identification in the exchange graphs, since certain established sets of permutations are not mutation equivalent to the initial seed.

Simple machine learning architectures could successfully learn to differentiate cluster algebras from their seeds, especially well between algebras of different types (finite, finite-mutation, infinite).
In these investigations we choose to use the cluster algebra seeds as inputs, as opposed to individual cluster variables, since there is such a large amount of repetition of individual variables between algebras.
However, further work may wish to examine the cluster variables as direct inputs also.

Other further work would naturally extend exchange graph structure analysis to different algebras of different ranks, and to higher depths where appropriate.
One may also be interested to analyse and ML the cluster complexes (dual to the generalised associahedra).
Finally, since the exchange graphs and quivers take a graphical form, graph neural networks may be a sensible more sophisticated architecture to  implement next.

\section*{Acknowledgement}
The authors wish to thank Gregg Musiker for invaluable insight and inspiring discussion.
PPD would like to thank the LMS for grants 42035 and 42111, and York St John University for grant QR21-22-63.
YHH would like to thank STFC for grant ST/J00037X/2.
E.~Heyes would like to thank SMCSE at City, University of London for the PhD studentship, as well as the Jersey Government for a postgraduate grant.
E.~Hirst would like to thank STFC for a PhD studentship.

\addcontentsline{toc}{section}{References}
\bibliographystyle{utphys}
\bibliography{references}

\end{document}